\documentstyle[12pt]{article}

\newtheorem{theo}{Theorem}[subsection]
\newtheorem{cor}[theo]{Corollary}
\newtheorem{prop}[theo]{Proposition}

 \newcommand{\beq}{\begin{equation}}
\newcommand{\eeq}{\end{equation}}
 \newcommand{\beth}{\begin{theo}}
\newcommand{\eth}{\end{theo}}

\makeatletter
\@addtoreset{equation}{section}
\makeatother

\newcommand{\Exmp}{{\it Example}}

\newcommand{\then}{\Rightarrow}
\def\Rn{{\bf R}^n}

\def\R{{\bf R}}
\def\C{{\bf C}}
\def\Cn{{\bf C}^n}
\def\Zn{{\bf Z}^n}
\def\Cm{{\bf C}^m}
\def\Cp{{\bf C}^p}
\def\Cq{{\bf C}^q}
\def\dpq{{\cal D}_{p,q}(\Omega)}
\def\dpqp{{\cal D}_{p,q}'(\Omega)}
\def\dpp{{\cal D}_{p,p}'(\Omega)}
\def\dpl{{\cal D}_p^+(\Omega)}
\def\Subset{\subset\subset}
\def\Luxt{\Lambda(u,x,t)}
\def\luxt{\lambda(u,x,t)}
\def\Luxa{\Lambda(u,x,a)}
\def\luxa{\lambda(u,x,a)}
\def\nux{\nu(u,x)}
\def\ntx{\nu(T,x)}
\def\ntp{\nu(T,\varphi)}
\def\nuxa{\nu(u,x,a)}
\def\okn{{1\le k \le n}}
\def\str{\sigma_T(x,r)}
\def\brv{{B_r(\varphi)}}
\def\srv{{S_r(\varphi)}}
\def\ntvr{\nu(T,\varphi,r)}
\def\ntv{\nu(T,\varphi)}
\def\vp{\varphi}
\def\mrv{\mu_r^\varphi}

\newcommand{\supp}{{\rm supp}\,}
\newcommand{\aoP}{\gamma^\varphi}
\newcommand{\calr}{{\cal R}}

\begin{document}

\title{ Singularities of plurisubharmonic functions
and positive closed currents}
\author{Alexander Rashkovskii}
\date{}
\maketitle{}

\begin{abstract}
This is a survey of results, both classical and recent, on
behaviour of plurisubharmonic functions near their
$-\infty$-points, together with the related topics for positive
closed currents.

\medskip

{\em Key  words:} plurisubharmonic function, positive closed
current, Lelong number, directional Lelong number, generalized
Lelong number, Monge-Amp\`ere operator

\end{abstract}

\bigskip

The history of plurisubharmonic functions starts in 1942 due to
P.~Lelong \cite{Lel} and K.~Oka \cite{O}. A model example of a
plurisubharmonic function is logarithm of modulus of a holomorphic
mapping $f$, and the behaviour of the mappings near their zeros
(or at least of $|f|$) corresponds to the behaviour of
plurisubharmonic functions near the points where they take the
value $-\infty$ (points of their singularities).

Being considered as elements of a distribution space, functions
$\log|f|$ serve as potentials for the zero sets of $f$. To this end,
the machinery of currents (due to de Rham) was developed by Lelong
for complex spaces, and central role here is played by closed
positive currents.
In 1957, Lelong proved that the trace measure of any
 closed positive current has density at every point of its support
\cite{Le0}. The main objects of his study were integration
currents over analytic varieties, and later R.~Thie showed that
for this case the densities coincide with the multiplicities of
the varieties (and he also called these values {\it Lelong
numbers}) \cite{Th}. The notion has turned out to be of great
importance. In particular, it provides us with a powerful link
between analytical and geometrical objects of modern complex
analysis. See Lelong's view of the subject in \cite{Le10},
\cite{keyLe11}. A collection of his relevant papers is presented
in \cite{Le12}.

Further developments in the field rest mainly on technique of
Monge-Amp\`ere operators, the key contribution being made by
J.P.~Demailly. Among various applications we mention those to
algebraic geometry and number theory (e.g., \cite{Bo}, \cite{D11},
\cite{LeGr}).


Here we present a  survey on the theory of Lelong numbers and
related topics (and not on the theory of plurisubharmonic
functions and positive closed currents). For more aspects of
pluripotential theory, see \cite{B}, \cite{C}, \cite{Kis5},
\cite{Za}, and for its
 backgrounds, \cite{Le}, \cite{Kl}, \cite{Ro0} and
lecture notes \cite{BlLN} and \cite{KoLN}. An excellent reference
is a book by Demailly \cite{Dbook}. The theory of positive closed
currents was surveyed in \cite{Sk4}; a systematic treatment of the
general theory of currents is contained in \cite{Fe}. The main
source for the present notes was the paper \cite{D1} (which is
actually a part of the book \cite{Dbook}). Closely connected to
the subject is the notion of pluricomplex Green functions; it is
left however beyond the scope of the survey, a brief presentation
being given in \cite{Kis5}.


Section \ref{sect:LN} introduces Lelong numbers of
plurisubharmonic functions as characteristics of their behaviour
at logarithmic singularity points. It leads to the notion of
Lelong number for positive closed currents (Section~2). In
Section~3, generalized Lelong numbers due to Demailly are studied.
Section~4 deals with analyticity theorems for Lelong numbers. In
Section~5, structural formulas for positive closed currents are
obtained. Finally, Lelong numbers for Monge-Amp\`ere currents are
evaluated in Section~6 by means of technique of local indicators.

Most of the results of Sections \ref{sect:LN} and 2 are classical
(see, e.g., \cite{Le}, \cite{LeGr}, \cite{GH}, \cite{Ho}), however
some recent developments (mainly due to C.O.~Kiselman) are
included as well. The exposition of Sections 3 -- 5 follows the
lines of \cite{D1}, with a few modifications. Section~6 is based
on  \cite{LeR}, \cite{R}, \cite{R4}.

\bigskip

{\small{\it Acknowledgments.} The notes are based on a course
given by the author at Summer Semester on Complex Analysis at Feza
Gursey Institute, Istanbul, 1999. The main part was written during
the author's stay at Mid Sweden University, Sundsvall, 2000. Many
thanks to the both institutions for their kind hospitality and
valuable discussions.}

\tableofcontents

\medskip

\section{Lelong numbers for plurisubharmonic functions}
\label{sect:LN}

The most transparent is the case of integration currents over analytic
 varieties $\{f(z)=0\}$ of codimension $1$,
and in fact it was studied by Lelong already in 1950 (\cite{Le2}).
The main idea is to reduce the problem to asymptotic behaviour of
the functions $\log|f|$. The same approach works for arbitrary
positive closed currents of bidegree $(1,1)$, and such currents
have local potentials -- plurisubharmonic functions. So we start
with the notion of Lelong number for plurisubharmonic functions.

\subsection{Plurisubharmonic functions}
\label{ssect:psh}

Throughout the exposition, $\Omega$ is a domain in $\Cn$, $n>1$,
and $u$ is a {\it plurisubharmonic (psh)} function in $\Omega$,
i.e., an upper semicontinuous function whose restriction to each
complex line $L$ is subharmonic in $\Omega\cap L$.

The class of all psh functions in $\Omega$ will be denoted by
$PSH(\Omega)$. Any psh function is locally integrable, and
the topology of $PSH(\Omega)$ is generated by
$L^1_{loc}$-convergence or, equivalently, by the weak convergence on
compactly supported continuous (or smooth) functions on $\Omega$.

Plurisubharmonicity of an upper semicontinuous function $u$ is characterized
by positive semidefineteness of its {\it Levi form}
$$
\sum_{j,k=1}^n \frac{\partial^2u}{\partial z_j\partial\bar z_k}
\eta_j\bar\eta_k\ge 0,\quad\forall \eta\in\Cn,
$$
the derivatives being understood in the sense of distributions.

For general theory of psh functions, see, e.g.,  \cite{Ho},
 \cite{Kl}, \cite{Le}, \cite{Ro0}.

Plurisubharmonicity can be viewed as convexity with respect to the
complex structure, however in contrast with convex functions on real affine
spaces, psh functions need not be continuous and, which is more important,
they main attain $-\infty$ values. And
 we will be interested mainly in the behaviour of psh functions $u$ near their
singularity points, that is, the points $x$ where $u(x)=-\infty$.
Note that for the functions $\log|f|$, such points form analytic varieties.

Some more terminology. A set $E\subset\Omega$ is called {\it
pluripolar} in $\Omega$ if there exists a psh function
$u\not\equiv -\infty$ in $\Omega$, such that $u|_E=-\infty$. It
was shown by B.~Josefson \cite{J} (see also \cite{BT2}) that if
$E$ is pluripolar in a neighbourhood of every its point then its
pluripolar in $\Omega$. A set $E\subset\Omega$ is called {\it
completely pluripolar} if there exists $u\in PSH(\Omega)$,
$u\not\equiv -\infty$, such that $u^{-1}(-\infty)=E$.

Most of our considerations are local, however sometimes we will have
to specify a domain $\Omega$ to be {\it pseudoconvex}. This means that there
exists a function $u\in PSH(\Omega)$ such that $u(z)\to+\infty$
as $z\to\partial\Omega$.

\subsection{Standard characteristics of local behaviour}

We  use the following notation:

$B_r(x)=\{x\in\Cn:|x|<r\}$, $S_r(x)=\partial B_r(x)$, $B_r:=B_r(0)$,
$S_r:=S_r(0)$;

$\tau_p=\pi^p/p!$ is the $2p$-volume of the unit ball in $\Cp$;

$\omega_p=2\pi^p/(p-1)!$ is the $(2p-1)$-volume of the unit sphere in $\Cp$.

\medskip
Let $u\in PSH(\Omega)$. For
 $x\in\Omega$ and $t<\log {\rm dist}\,(x,\partial\Omega)$,
define
$$
\Lambda(u,x,t)=\sup\,\{u(z):z\in B_{e^t}(x)\},
$$
$$
\lambda(u,x,t)=\omega_n^{-1}\int_{S_1}u(x+ze^t)\,dS_1(z).
$$
With $t$ fixed, the both functions are continuous and plurisubharmonic
in $x$.
With $x$ fixed, they are convex and increasing in $t$,
$$
u(x)\le \luxt\le\Luxt,
$$
 and
$\luxt,\,\Luxt\to u(x)$ as $t\to -\infty$. Moreover,
$\luxt/\Luxt\to 1$ as $t\to -\infty$ if $u(x)=-\infty$
(a consequence of Harnack's inequality). For more relations between
these functions, see Sections \ref{ssect:1.6} and \ref{ssect:1.7}.

\subsection{Definition of Lelong number for
plurisubharmonic functions}

Since psh functions are locally integrable, it is possible to apply
the machinery of differential operators.

Let $\Delta=4\sum_k\partial^2/\partial z_k\partial\bar z_k$
be the Laplace operator,
then $\Delta u$ is a positive
measure on $\Omega$ (which is, up to a constant factor,
the {\it Riesz measure} of $u$ considered as a subharmonic function in
${\bf R}^{2n}$).
 Denote
$$
\sigma_u(x,r)={1\over 2\pi}\int_{B_r(x)}\Delta u.
$$
Green's formula implies

\begin{prop}
\label{prop:deriv}
\beq
\label{eq:deriv}
{\sigma_u(x,r)\over
\tau_{n-1}r^{2n-2}}=\frac{\partial\lambda(u,x,\log r)}{\partial\log
r},
\eeq
$\partial/\partial\log r$ being understood as the left derivative.
\end{prop}
\medskip

Since $\luxt$ is convex and increasing, the right-hand side of
(\ref{eq:deriv}) is increasing in $r$, so is its left-hand side and hence
there exists the limit
\beq
\label{eq:LN}
        \lim_{r\to 0} {\sigma_u(x,r)\over
\tau_{n-1}r^{2n-2}}=:\nux,
\eeq
{\it the Lelong number of $u$ at $x$}.
In other words, the Lelong number of $u$ is the $(2n-2)$-dimensional
density of its Riesz measure at $x$.

This point of view on Lelong numbers will be developed in Section~2.
And here we concentrate on representations of $\nux$ in terms
of the asymptotic behaviour of $u$ near $x$. Consideration of the
right-hand side of
(\ref{eq:deriv}) and the equivalence between $\luxt$ and $\Luxt$ gives us

\beth \label{theo:nu} {\rm \cite{Av}, \cite{Kis1}} \beq
\label{eq:nu} \nux = \lim_{t\to -\infty}{\luxt\over t}= \lim_{t\to
-\infty}{\Luxt\over t}. \eeq \eth

As follows from this result, $\nux>0$ is possible only when
$u(x)=-\infty$, and the converse is not true. Namely,
if the Lelong number of a psh function $u$ is strictly positive,
then $u$ has a logarithmic singularity at $x$:

\begin{cor}
\label{cor:lognu}
$$
\nux=\sup\,\{\nu>0:u(z)\le \nu\log|z-x|+O(1),\ z\to x\}.
$$
\end{cor}

Evidently, for finite collections of psh functions $u_k$,
$$
\nu(\sum_k u_k,x)=\sum_k\nu(u_k,x)
$$
 and
$$
\nu(\max_k\,u_k,x)=\min_k\,\nu(u_k,x).
$$

Next result is more difficult, and it is a particular case of
Theorem~\ref{theo:pback} or \ref{theo:invcur}.

\beth \label{theo:inv} {\rm \cite{Siu}}
 The Lelong number $\nux$ is independent of the
choice of local coordinates.
\eth

\subsection{Examples}

The following relations can be easily derived from Theorem~\ref{theo:nu}.

(a) if $u(z)=\log|z|$, then $\nu(u,0)=1$.

(b) let $u(z)=\log|f(z)|$ and $f:\Omega\to{\bf C}$ be a
holomorphic function, $f(x)=0$; then $\nu(u,x)=m$, the
multiplicity (vanishing order) of $f$ at $x$ (the least degree of
a monomial in the Taylor expansion of $f$ near $x$).

(c) if $u(z)=\log|f(z)|={1\over 2}\log\sum_k|f_k|^2$ and
$f=(f_1,\ldots,f_m):\Omega\to\Cm$ is a holomorphic mapping,
$f(x)=0$, then
$$
\nux=\min_k\,m_k,
$$
where $m_k$ are the multiplicities of the
zeros of the components $f_k$ of the mapping $f$ at $x$.

\subsection{Lelong numbers of slices and pull-backs}

Fix $x\in\Omega$. Given $y\in\Cn$, $L$ is the complex line
through $x$ and $y$, and $u_y$ is the restriction of $u$ to
$\Omega\cap L$ (the slice of $u$ on $L$):
\beq
\label{eq:defslice}
u_y(\zeta):=u(x+\zeta y)
\in SH(\Omega\cap L), \quad\zeta\in{\bf C}.
\eeq

\beth
\label{theo:slice0}
 $\nu(u_y,0)\ge\nux$ for all $y\in\Cn$, and
$\nu(u_y,0)=\nux\ \forall y\in\Cn\setminus A$, $A$ being a pluripolar
subset of $\Cn$.
\eth

The first statement is evident in view of (\ref{eq:nu}), and the second
can be derived then from the relation
\beq
\label{eq:Croft0}
\nux=\omega_n^{-1}\int_{S_1} \nu(u_y,0)\,dS_1(y).
\eeq

On the exceptional set $A$, the values $\nu(u_y,0)$ can behave as bad
as possible:

\begin{prop}
{\rm \cite{CPo}} For any countable $G_\delta$-subset $\{L_j\}$ of
the Riemann sphere and any sequence $c_j>c>0$ there exists a psh
function $u$ in the unit ball of ${\bf C}^2$ such that
$\nu(u_{y_j},0)=c_j$ and $\nu(u,0)=c$.
\end{prop}

Theorem \ref{theo:slice0} remains true when considering slicing of
the functions by $p$-dimensional planes \cite{Kis9}.

Let now $f$ be a holomorphic mapping $\Omega'\to\Omega$ with
$f(x')=x$, and $f^*u$ be
the {\it pull-back} of a function $u\in PSH(\Omega)$, that is, $f^*u(z)=
u(f(z))$.

\beth \label{theo:pback} {\rm\cite{Kis9}} $\nu(f^*u,x')\ge \nux$.
Moreover, strict inequality here occurs only for $f$ belonging to
a polar set in the Frechet space of holomorphic mappings on a
neighbourhood of $x'$. \eth

A relation in the opposite direction
is given by

\beth \label{theo:favr} {\rm \cite{F}, \cite{Kis6}} If $f(U)$ has
non-empty interior for every neighbourhood $U$ of $x$, then there
exists a constant $C$, independent of $u$, such that
 $\nu(f^*u,x')\le C\nux$ for any function $u$ plurisubharmonic in
a neighbourhood of $x$. No such bound is possible if $f(U)$
has no interior points for some neighbourhood $U$.
\eth

The proof can be derived, for example, from
Proposition \ref{prop:sublevel} below.

Characteristics for singularities of pull-backs of psh functions
under blow-ups of points and subvarieties of $\Omega$, {\it
microlocal Lelong numbers}, were defined and studied in
\cite{Abr}, see also \cite{Lau}.

\subsection{Attenuating the singularities}
\label{ssect:1.6}

The following construction was proposed by C.O.~Kiselman
\cite{Kis1}, \cite{Kis3}. Let $u,q\in PSH(\Omega)$, and $q(x)\ge
-\log{\rm dist}\,(x,\partial\Omega)$. For any $x\in\Omega$ and
$\alpha>0$, define
$$
u_{\alpha,q}(x)=\inf\{\luxt-\alpha \, t:\: t<-q(x)\}.
$$
It turns out to be a plurisubharmonic function in $x$ (by the
minimum principle of Kiselman \cite{Kis8}) and concave in
$\alpha$. Moreover, its singularities are close to those of $u$:

\begin{theo}
\label{prop:atten} {\rm\cite{Kis1}, \cite{Kis3}} If $\nu(q,x)=0$
then $\nu(u_{\alpha,q},x)=\max\{\nux-\alpha,0\}$.
\end{theo}

The same is true when using the function $\Luxt$ instead of $\luxt$, however
the corresponding function $u_{\alpha,q}$ may be different. Generally, the
functions $\Luxt$ and $\luxt$ can be related as follows.

\begin{prop}
{\rm \cite{Le6}} \beq \label{eq:Ll}
\limsup_{t\to-\infty}\,[\Luxt-\luxt]\le c_n\,\nux \eeq with
$c_1=0$ and $c_n=\sum_{1\le k\le 2n-2}k^{-1}$ for $n>1$.
\end{prop}

Moreover, if $n>1$ and $\nux>0$, the relation
$\lim_{t\to-\infty}[\Luxt-\luxt]=0$
 is not true in general; this limit may even not exist.

\subsection{Principal parts}
\label{ssect:1.7}

Given a function $u\in PSH(\Omega)$, let $u_y$ be its slice (\ref{eq:defslice})
along the line through $x\in\Omega$ and $y\in\Cn$.
Consider
$$
\tilde u_x(y)=\limsup_{\zeta\to 0}[u_y(\zeta)-\nux\,\log|\zeta|].
$$
The function $u$ is said to have a  principal  part at $x$ if
$\tilde u_x\not\equiv-\infty$, and the function $\tilde u_x$ is
called the {\it  principal  part} of $u$ at $x$ \cite{Le5}.

\begin{theo}
{\rm\cite{Le5}} A function $u$ has a principal part at $x$ if and
only if and only if there exists
$a_{u,x}=\lim_{t\to-\infty}[\Luxt-\nux\,t]>-\infty$.
\end{theo}

Note that if $u$ has a principal part then there exist the limit in
(\ref{eq:Ll}). Besides, the function $\tilde u_0(y) - a_{u,0}$
is the maximal element of the limit set $LS(u)$ of $u$ as defined in
Section~\ref{ssect:tangents}.

\subsection{Directional Lelong numbers}

The preceding considerations are based on comparing psh functions
with convex functions of a real variable. More detailed information
can be obtained by comparing them with convex functions in $\Rn$.

For $x\in\Omega,\ a=(a_1,\ldots,a_n)\in\Rn$, one can consider the
polydisk characteristics
$$
\luxa:=(2\pi)^{-n}\int_{[0,2\pi]^n} u(x_k+e^{a_k+i\theta_k})\,
d\theta,
$$
$$
\Luxa:=\sup\,\{u(z):z\in T_a(x)\},
$$
where
$$
T_a(x)=\{z:|z_k-x_k|=e^{a_k},\:1\le k\le n\}.
$$

These functions are convex in $a$ and increasing in each $a_k$,
$u(x)\le \luxa\le\Luxa$,
 and
$\luxa,\,\Luxa\to u(x)$ as $a_k\to -\infty$, $1\le k\le n$.

So, there exist the limits \beq \label{eq:dir} \lim_{t\to
-\infty}{\lambda(u,x,ta)\over t}= \lim_{t\to
-\infty}{\Lambda(u,x,ta)\over t}=:\nuxa \eeq for any $
a\in{\R}_+^n$, and the value $\nuxa$ is called the {\it
directional} (or {\it refined}) {\it Lelong number} due to
Kiselman \cite{Kis2}, \cite{Kis3}.

If ${\bf 1}=(1,\ldots,1)$ then $\nux=\nu(u,x,{\bf 1})$; besides,
\beq
\label{eq:dircomp}
{\min_j a_j}\,\nux\le\nuxa\le{\max_j a_j}\,\nux.
\eeq

\medskip
\Exmp. Let
$u(z)=\log|f(z)|$ with a holomorphic function $f$, $f(x)=0$.
In a neighborhood of $x$ the function has the form
$$
f(z)=\sum_{J\in\omega_x} c_J(z-x)^J, \quad c_J\neq 0
$$
($\omega_x\subset\Zn_+$). Then (\cite{Le7}) \beq \label{eq:index}
\nuxa=\min\,\{\langle a,J\rangle:\: J\in\omega_x\}. \eeq Note that
when $f$ is a polynomial, the right-hand side of (\ref{eq:index})
is known as the {\it index} of the polynomial with respect to the
weight $a$ \cite{La}.

\medskip

Directional Lelong numbers appear naturally when studying monomial
transforms $f:z\mapsto (z^{M_1},\ldots,z^{M_m})$, $M_j$ being the
$j$-th row of a matrix $M$ with nonnegative integer entries, $\det
M\neq 0$: $ \nu(f^*u,0)=\nu(u,0,{\bf 1}M^*)$. More generally, for
any direction $a$, we have $ \nu(f^*u,0,a)=\nu(u,0,aM^*)$, so the
monomial transforms affect the directional numbers linearly
\cite{Kis3}.

The procedure of attenuating singularities (Section~\ref{ssect:1.6})
can be applied as well
to the directional Lelong numbers by considering the function
$$
u_{a,\alpha,q}(x)=\inf_{t<-q(x)} (\Lambda(u,x, ta)-\alpha\,t).
$$

\begin{prop}
\label{prop:attendir} {\rm\cite{Kis3}} If $\nu(q,x)=0$ then
$\nu(u_{a,\alpha,q},x,a)=\max\{\nuxa-\alpha,0\}$.
\end{prop}

\subsection{Partial Lelong numbers}\label{ssec:pln}

When taking the mean values with respect to some of the variables
$z_k$, it produces the {\it partial Lelong numbers} \cite{Le3},
\cite{Kis3}, \cite{W}. Namely, for an ordered $p$-tiple
$J=(j_1,\ldots,j_p)$ from $\{1,\ldots,n\}$ and a direction
$a'\in\R_+^p$, \beq \label{eq:part} \nu_J(u,x,a'):=\lim_{t\to
-\infty}t^{-1}\lambda_J(u,x,ta'),\quad a'\in\R_+^p, \eeq
$\lambda_J(u,x,a')$ being the mean value of $u$ over the set
$$
T_{a'}(x)=
\{z\in\Cn:\: |z_{j_m}-x_{j_m}|=\exp{a'_{j_m}}\:
(1\le m\le p),\ z_k=x_k\: (k\not\in J) \},\quad
a'\in\R^p.
$$

It was shown in \cite{W} that if $u$ satisfies certain regularity
conditions near $x$, then $\nuxa\to\nu_J(u,x,a')$ as $a_k\to
+\infty$, $k\not\in J$, and $a_j=a'_j$, $j\in J$. In the general
situation the limit exists although it can be strictly less than
the corresponding partial Lelong number, even if $u$ is assumed to
be locally bounded outside $x$ and multicircled around $x$ (i.e.,
$u(x+\zeta)=u(|x_1+\zeta_1|,\ldots,|x_n+\zeta_n|)$ for any $\zeta$
near $0$). Indeed, for
$$u(z_1,z_2)=\max\,\{\log|z_1|,-|\log|z_2||^{1/2}\},$$
$\nu(u,0,a)=0$ for all $a\in\Rn_+$, while $\nu_1(u,0,a')=a'$.

\subsection{Sublevel sets, integrability index, and multiplier ideals}
\label{ssect:1.9}

The rate of approaching $u(z)$ its $-\infty$ value can
be characterized by behaviour of the volumes of its
{\it sublevel sets}
$$
A_{u}(t)=\{y: u(y)<t\}
$$
as $t\to-\infty$.

Another way is to study local integrability of $e^{-u/\gamma}$
for $\gamma>0$ (note
that since $u$ may be discontinuous, even integrability of $e^{-u}$
near $x$ with $u(x)>-\infty$ is far from being evident).
The value
\beq
\label{eq:intind}
I(u,x)=\inf\,\{\gamma>0:e^{-u/\gamma} \in L^2_{loc}(x)\}
\eeq
is called the {\it integrability index}, or {\it Arnold multiplicity},
of $u$ at $x$.

The both characteristics turn out to be closely related.

\begin{theo}
\label{prop:sublevel} {\rm \cite{DK}, \cite{Kis6}} $I(u,x)$ is the
infimum of $\gamma>0$ such that $e^{-2t/\gamma}{\rm
Vol}\,A_u(t)\cap U$ is bounded as $t\to -\infty$ for some
neighbourhood $U$ of the point $x$.
\end{theo}

In other words, volumes of the sublevel sets have exponential
decay with the rate controlled by the integrability index. Uniform
relations of such type for families of psh functions were obtained
in \cite{Zer}. Besides, Proposition~\ref{prop:CLN} implies that
the size of the sublevel sets can be also measured in terms of
plurisubharmonic capacity (\ref{eq:capacity}):

\begin{prop}
Given sets $U\Subset V\Subset\Omega$, there exists a constant $C$ such
that for any negative function $u\in PSH(\Omega)$ and $t<0$,
$$
Cap\,(A_u(t)\cap U,\Omega)\le \frac{C}{|t|}\|u\|_{L^1(V)}.
$$
\end{prop}

Integrability indices can be estimated in terms of Lelong numbers
\cite{Sk}:
$$
{1\over n}\nux\le I(u,x)\le \nux,
$$
the extremal situation being realized for $u=\log|z_1|$ and $u=\log|z|$.
A more refined relation is given by

\beth \label{theo:intdir} {\rm \cite{Kis3}} $\sup\,\{\nuxa:
\sum_ja_j=1\}\le I(u,x)\le\nux$. Moreover, if
$u(z)=u(|z_1-x_1|,\ldots,|z_n-x_n|)$ near $x$, then
$I(u,x)=\sup\,\{\nuxa: \sum_ja_j=1\}$.
\end{theo}

An analytic object representing the singularities of $u$ is the
{\it multiplier ideal sheaf} ${\cal J}(u)$ consisting of germs of
holomorphic functions $f$ such that $|f|e^{-u}\in L^2_{loc}$. It
is a coherent analytic sheaf \cite{Nad}. Moreover, if
$U\Subset\Omega$ is pseudoconvex, then the restriction of ${\cal
J}(u)$ to $U$ is generated as an ${\cal O}_U$-module by a Hilbert
basis $\{\sigma_l\}$ of the Hilbert space $H_u(U)$ of holomorphic
functions $f$ on $U$ such that $|f|e^{-u}\in L^2(U)$ \cite{DeLL}.

An application of this notion will be given in Section~\ref{ssect:4.4}.


\section{Lelong numbers for positive closed currents}
\label{sect:currents}

Up to this moment we developed an approach to Lelong numbers based
on asymptotic properties of plurisubharmonic functions. More information
can be obtained by considering them as densities of the
Riesz measures. To do this, the measures should be viewed as the trace
measures of the corresponding positive closed
currents of degree $(1,1)$. And this can be
extended to currents of higher degrees.  One of the motivations for such
an extension is as follows. When $u=\log|f|$ and $f:\Omega\to\C$,
the Lelong number of $u$ at a point $x$ is just
the multiplicity of the zero of $f$ at $x$, and it is not the case
for holomorphic mappings $f$ to $\Cm$ when $m>1$. As will be seen
in Section~5, multiplicities of holomorphic mappings
can be characterized as Lelong numbers of certain currents of higher degrees.

So we pass to Lelong numbers for positive closed currents,
starting with recalling some basic notions of the theory of
currents. The subject of Sections 2.1 -- 2.5 is treated, e.g., in
\cite{Ho},  \cite{Kl}, \cite{Le}, \cite{LeGr}. More information on
Lelong numbers of the currents will be presented in Sections 4 and
5.

\subsection{Positive closed currents}

Let $\dpq$ be the space of all smooth compactly supported differential
forms $\phi$ of bidegree $(p,q)$ on $\Omega$:
$$
\phi=\sum_{|I|=p,|J|=q}\phi_{IJ}dz_I\wedge d\bar z_J,\quad \phi_{IJ}\in
{\cal D}(\Omega),
$$
with the topology of $C^\infty$-convergence.

Differentiation of forms: the operator $\partial: \dpq\to{\cal
D}_{p+1,q}(\Omega) $ is defined by
$$
\partial\phi=\sum_{I,J}\sum_{1\le k\le n}
{\partial\phi_{IJ}\over \partial z_k}dz_k\wedge dz_I\wedge d\bar z_J,
$$
and
$\bar\partial: \dpq\to{\cal D}_{p,q+1}(\Omega) $ is given by
$$
\bar\partial\phi=\sum_{I,J}\sum_{1\le k\le n}
{\partial\phi_{IJ}\over \partial \bar z_k}d\bar z_k\wedge dz_I\wedge
d\bar z_J.
$$

The {\em currents} of bidimension $(p,q)$ (bidegree $(n-p,n-q)$) are
elements of the dual space
$\dpqp$ (continuous linear functionals on $\dpq$).
Any current $T\in\dpqp$ has a representation
$$
T=\sum_{|I|=n-p,|J|=n-q}T_{IJ}dz^I\wedge d\bar z^J,\quad T_{IJ}\in
{\cal D}'(\Omega).
$$
The action of $T$ on $\phi$ will be written as
$\langle T,\phi\rangle$ or $\int T\wedge\phi$.

The topology on $\dpqp$ (referred to as {\it the weak topology
of currents}):
$$
T_j\to T\iff  \langle T_j,\phi\rangle
\to \langle T,\phi\rangle\quad\forall\phi\in\dpq.
$$

Differentiation of currents:
$$
\langle \partial T,\phi\rangle:=(-1)^{p+q+1}\langle
T,\partial\phi\rangle,\quad
\langle \bar\partial T,\phi\rangle:=(-1)^{p+q+1}\langle
T,\bar\partial\phi\rangle.
$$
The operators
$$
d=\partial+\bar\partial,\quad
d^c={\partial-\bar\partial\over 2\pi i}
$$
are real, and $dd^c={i\over\pi}\partial\bar\partial$. (There is no
general convention on normalizing the operator $d^c$, some authors
use $d^c=i(\bar\partial-\partial)$; we prefer the above one to
avoid extra factors $(2\pi)^p$ in the sequel.)



A current $T\in\dpp$ is called {\it positive} ($T\ge 0$) if
$\langle T,\phi\rangle\ge 0$ for every differential form
$\phi=i\alpha_1\wedge\bar\alpha_1\wedge\ldots\wedge
i\alpha_p\wedge\bar\alpha_p$ with $\alpha_k\in{\cal
D}_{1,0}(\Omega)$. (In the literature, such currents are sometimes
called {\it weakly positive}).
 The coefficients $T_{IJ}$ of such a current $T$ are Borel
measures on $\Omega$. Therefore, the action of a positive current
$T=\sum T_{IJ}dz^I\wedge d\bar z^J$ can be continuously extended
to the space of compactly supported forms $\phi$ with {\sl
continuous} coefficients $\phi_{KL}$ and
$$|\langle T,\phi\rangle|\le \|T\|_{\supp\phi}\|\phi\|,$$
where $\|T\|_E=\sum |T_{JK}|_E$, $|T_{JK}|_E$ is the total
variation of the measure $T_{JK}$ on $E$ and
$\|\phi\|=\sup_{K,L,x}|\phi_{KL}(x)|$.

Denote by
$$
\beta:={i\over 2}\sum_\okn dz_k\wedge d\bar z_k={\pi\over
2}dd^c|z|^2
$$
the standard K\"ahler form on $\Cn$, so
$$
\beta_p:={1\over p!}\beta^p
$$
is the $p$-dimensional volume element. Then for every positive
current $T\in\dpp$,
$$\|T\|_E\le c_n |T\wedge\beta_p|_E.$$

A current $T$ is called {\it closed} if $dT=0$. When $T\in\dpp$,
this is equivalent to saying that $\partial T=0$ or $\bar\partial
T=0$.

$\dpl$ will denote the cone of all positive closed currents from
$\dpp$.

An important tool in the theory of positive closed currents is the
following Skoda-El Mir extension theorem.

\beth \label{theo:SkEM} {\rm \cite{Sk3}, \cite{EM}, \cite{Si}} Let
$E$ be a closed complete pluripolar set in $\Omega$ and $T\in{\cal
D}_p^+(\Omega\setminus E)$ whose coefficients $T_{IJ}$ have
locally finite mass near $E$. Consider the current $\tilde
T=\sum\tilde T_{IJ}dz^I\wedge d\bar z^J$ with the measures $\tilde
T_{IJ}(A):=T_{IJ}(A\setminus E)$ for all Borel $A\subset\Omega$.
Then $\tilde T\in\dpl$. \eth

(The current $\tilde T$ is called the {\it simple}, or {\it trivial},
{\it extension} of $T$, and it was first introduced by Lelong when
studying integration over analytic varieties, see Example~3 in
Section~\ref{ssect:2.2}).

\medskip

If $f:\Omega\to\Omega'\subset\Cn$ is a holomorphic mapping
such that its restriction
to the support of a current $t\in\dpl$ is proper, then the {\em direct image}
(or {\em push-forward}) $f_*T$  of $T$ is defined by the relation
$\langle f_*T,\phi\rangle=\langle T,f^*\phi\rangle$. If $T\in\dpl$ then
$f_*T\in {\cal D}_p^+(\Omega')$.

If a holomorphic mapping $f:\Omega\to\Omega'$ has constant maximal rank
on $\Omega$, then for any current $T\in {\cal D}'_{p,l} (\Omega')$
its {\em inverse
image} (or {\em pull-back}) $f^*T$ is defined as
$\langle f^*T,\phi\rangle=\langle T,f_*\phi\rangle$.

The inverse image of $T$ can be also defined for any surjective
holomorphic mapping $f$ with $\dim f^{-1}(x)=0$ for every
$x\in\Omega'$ (see \cite{Meo}).


\subsection{Examples of currents}
\label{ssect:2.2}

The standard examples are as follows.

1) Currents generated by psh functions:
$$
u\in PSH(\Omega)\iff \left(\frac{\partial^2u}{\partial
z_j\partial\bar z_k}\right)\ge 0 \iff dd^cu\in{\cal
D}_{n-1}^+(\Omega).
$$
Furthermore, if $T\in{\cal D}_{n-1}^+(\Omega)$ then for any
$x\in\Omega$ there is a neighbourhood $U_x$ and a function $u_x\in
PSH(U_x)$ such that $T=dd^cu_x$ in $U_x$.

2) For $M$ a complex manifold of dimension $p$, the current $[M]$
of integration over $M$ is defined as
$$
\langle [M],\phi\rangle=\int_M\phi.
$$
Then $[M]\in\dpl$ (that it is closed, follows from Stokes' theorem).

3) Integration currents over analytic varieties. Let
$A$ be an analytic variety, i.e., locally
$A=\{z: f_\alpha(z)=0,\ \alpha\in{\cal A}\}$,
and $Reg\, A$ be the set
of its regular points (where $A$ is locally a manifold).
If $A$ is of pure dimension $p$, define
$$
\langle [A],\phi\rangle:=\int_{Reg\,A} \phi.
$$
Then $[A]\in\dpl$. (Non-trivial part is that $[A]$ is closed; this
fundamental result is due to P.~Lelong \cite{Le0}, and it can be
seen today as a consequence of Theorem~\ref{theo:SkEM}.)

4) Holomorphic chains $T=\sum \alpha_k [A_k]\in\dpl$, where
$\alpha_k\in {\bf Z}_+$ and $A_k$ are analytic varieties of pure
dimension $p$. When $p=n-1$, the holomorphic chains represent
positive, or effective, divisors.

\subsection{Lelong numbers for currents}

For $T\in\dpl,\ \sigma_T:=T\wedge\beta_p\in{\cal D}_0^+$ is the
{\it trace measure} of $T$.
(If $T=dd^cu$ then $\sigma_T$ is just the Riesz measure of $u$.)

Denote $\str=\sigma_T(B_r(x))$. It can be also
represented in the following form.

\begin{prop}
$$
\str=\tau_pr^{2p}\int_{B_r(x)}T\wedge\left(dd^c\log|z-x|\right)^p.
$$
\end{prop}

Therefore,
$$
\nu(T,x,r):={\str\over\tau_pr^{2p}}\searrow\ntx,
$$
{\it the Lelong number} of the current $T\in\dpl$ at $x$.
So,
\beq
\ntx=\lim_{r\to 0}{1\over\tau_pr^{2p}}\int_{B_r(x)}T\wedge
\beta_p=
\lim_{r\to 0}\int_{B_r(x)}T\wedge
\left(dd^c\log|z-x|\right)^p.
\eeq

The Lelong number of a current $T\in\dpl$ can be viewed as the
$2p$-dimensional density of its trace measure $\sigma_T$ or,
equivalently, as the mass charged at $x$ by its ``projective''
trace measure $T\wedge (dd^c\log|\cdot-x|)^p$.

\begin{cor}
\label{cor:volest}
 $\str\ge \tau_pr^{2p}\ntx$.
\end{cor}

\subsection{Special cases}

Lelong numbers of the model examples of currents presented
in Section~\ref{ssect:2.2} are as follows.

1) If $T=dd^cu$ with a plurisubharmonic function $u$, then $\ntx=\nux$.
This follows easily from the original definition (\ref{eq:LN}) of
Lelong numbers
of psh functions, since $\sigma_u=\sigma_T$. For the uniformity, it should
be written $\nu(dd^cu,x)$ instead of $\nux$, however we prefer to
keep the original notation for the Lelong numbers of functions, both
for the sake of brevity and since this is their standard notation.

2) The situation with manifolds is also easy: For any complex manifold
$M$, $\nu([M],x)=1$ at all points
$x\in M$ (and of course  $\nu([M],x)=0$ for $x$ outside $M$).

3) Much more difficult is the result about the Lelong numbers of
analytic varieties (Thie's theorem). Let $A\subset\Omega$ be an
analytic variety of pure dimension $p$, then $[A]\in\dpl$. For any
fixed $x\in A$ there is a neighbourhood $U$ and local coordinates
$(z',z'')\in{\bf C}^p\times{\bf C}^{n-p}$ such that $A\cap
U\subset\{(z',z''):\:|z''|\le C|z'|\},\ C>0$. In these
coordinates, $x=0$. Consider the projection $\pi: A\cap(U'\times
U'')\to U'$; it gives a ramified covering of $U'$. The number of
sheets of the covering, $m_x$, is called the {\it multiplicity} of
$A$ at $x$.

\beth \label{theo:thie} {\rm \cite{Th}} For an analytic variety
$A$ of dimension $p$, $\nu([A],x)$ equals the multiplicity $m_x$
of $A$ at $x$. \eth

For any Borel $D\Subset\Omega$,
the $2p$-dimensional volume of $A\cap D$ is precisely
$\sigma_{[A]}(D)$. Therefore, we have

\begin{cor}
{\rm (Volume estimation)} If $K\Subset A$ and $r_0<{\rm
dist}\,(K,\partial \Omega)$, then
$$
\tau_pr^{2p}m_x\le {\rm Vol}_{2p}\,A\cap B_r(x) \le C(r_0,K,A)\,r^{2p}
\quad\forall r<r_0,\ \forall x\in K.
$$
\end{cor}

\subsection{Stability of Lelong numbers}

Lelong numbers of currents, just as of psh functions, have the following
invariance property.

\beth \label{theo:invcur} {\rm \cite{Siu}} The Lelong number
$\ntx$ is independent of the choice of local coordinates. \eth

Another important relation concerns Lelong numbers of slices of
positive closed currents.
Given a current $T\in\dpl$, one can define its slices $T|_{x+S}$ by complex
planes $x+S$, where
 $x\in\Omega$ and $S\in G(q,n)$, $q\ge n-p$,
$G(q,n)$ being the Grassmanian of $q$-dimensional complex linear
subspaces of $\Cn$. For any $x\in\Omega$ and almost all (with
respect to the Haar measure) $ S\in G(q,n)$, the current
$T|_{x+S}:= T\wedge [x+S]\in{\cal D}_{p+q-n}^+(\Omega)$ is well
defined and supported on $x+S$. For details of the slicing theory,
see \cite{Fe}.

\beth \label{theo:slice} {\rm \cite{Siu}} Let $T\in\dpl$ and $q\ge
n-p$. Then for every $x\in\supp T$ and almost all $S\in G(q,n)$,
we have $\ntx=\nu(T|_{x+S},x)$. \eth

The proof follows along the same lines as that of Theorem~\ref{theo:slice0},
however technically is much more complicated. In particular, the role
of (\ref{eq:Croft0}) is played now by the Crofton formula
$$
\int_{G(p,n)}[S]\,d\kappa(S)=(dd^c\log|z|)^{n-p}
$$
($d\kappa$ is the Haar measure on $G(p,n)$).

\subsection{Tangents to positive closed currents}
\label{ssect:tangents}

Let $h_r(z)=z/r$, $r>0$. Given a current $T\in\dpl$, consider its
push-forwards $h_{r*}T$ (we suppose that $0\in\Omega$).
If $T$ is the integration current over an analytic
variety, then the limit
$$
\lim_{r\to 0}h_{r*}T=:tc(T)
$$
exists and is called the {\it tangent cone} to $T$ at $0$, see
\cite{Ha}. In particular, when $T=dd^c\log|f(z)|$ with $f$ a
holomorphic function near the origin, $tc(T) =dd^c\log|P_m|$,
where $P_m$ is the homogeneous term of $f$ of minimal degree.

A problem of existence of the tangent cones for other positive closed
currents was formulated by Harvey. It was shown by Kiselman
to be answered
in negative even in the situation of currents of bidegree
$(1,1)$. Given a psh function $u$ on
a neighbourhood of the origin, consider the family
\beq
\label{eq:tang}
u_r(z):=h_{r*}u(z)-\sup_{B_r}u=u(rz)-\sup_{B_r}u.
\eeq
It is relatively compact in $L_{loc}^1(\Cn)$
for $r<r_0$. Denote by $LS(u)$ the limit set
of this family as $r\to 0$, i.e.,
the set of all partial limits of sequences $u_{r_j}$, $r_j\to 0$.
Evidently, the current
$dd^c u$ has the tangent cone if and only if $LS(u)$ consists of a unique
element. For example,  $LS(\log|f|)=\log|P_m|-\sup_BP_m$.

It can be checked that $LS(u)$ is formed by functions $g\in PSH(\Cn)$
satisfying $g(tz)=\nu(u,0)\log|t|+g(z)$ for any $t\in\C$ (in particular,
$\nu(g,0)=\nu(u,0)$  $\forall g\in LS(u)$).
When $n=1$, $LS(u)$ consists of the only function
$\nu(u,0)\log|z|$ for every
$u$ (pluri)subharmonic near the origin of $\C$,
and it is not the case in higher dimensions:

\beth \label{theo:tangKis} {\rm \cite{Kis4}} Given any closed and
connected subset $M$ of $PSH(\Cn)$, $n>1$, which consists of
functions $g$ satisfying $\sup_B g<0$ and $g(tz)=C\log|t|+g(z)$
for all $t\in\C$, there exists a psh function $u$ with $LS(u)=M$
(and, consequently, $\nu(u,0)=C$). \eth

The corresponding limit set $LS(T)$ for $T\in\dpl$, $p\le n-1$,
 is formed by currents
$S$ with $\nu(S,0)=\nu(T,0)$ and $S\wedge (dd^c\log|z|)^p=0$ on
$\Cn\setminus\{0\}$ \cite{D10}. Theorem \ref{theo:tangKis} was
extended to currents of arbitrary degree in \cite{Blel}.

It turns out that the existence of a tangent cone depends on the rate
at which the projective masses $\nu(T,0,r)$ approach the Lelong number of $T$:

\beth {\rm\cite{BDM}} Let $T\in\dpl$. If the function
$n_T(r):=r^{-1}[\nu(T,0,r)-\nu(T,0)]$
is integrable at $0$, then $T$ has tangent cone at the origin.
\eth

Moreover, the condition is sharp: it is possible to construct a
current $T$ of bidegree $(1,1)$ with no tangent cone and such that
the integral of the function $n_T$ has the divergence rate at $0$
as small as one likes \cite{BDM}.

The set of points where a current $T\in{\cal D}_{n-1}^+$ has no
tangent cones was studied in \cite{Blel2}.

\subsection{Pluripositive currents and other generalizations}

A real current $T\in\dpp$ is called {\it pluripositive} if
$dd^cT\ge 0$. It was shown in \cite{Sk3} that the notion of Lelong
number can be extended to all negative pluripositive currents. The
idea was developed in \cite{Si}; in particular, it was applied to
positive pluripositive currents.  For other types of Lelong
numbers for such currents, see \cite{AB}. Extensions of
pluripositive currents across sets of small Hausdorff dimension
(analogs of \ref{theo:SkEM}) were studied in \cite{Si},
\cite{DaE}.

Being a local characteristic, the notion of Lelong number extends
easily to positive closed currents on complex manifolds. For
currents on singular complex spaces, see \cite{D3}. Currents on
almost complex manifolds were treated in \cite{Hag}.


\section{Generalized Lelong numbers due to Demailly}

An important notion of  {\it generalized Lelong numbers with
respect to psh weights} was introduced and studied by
J.P.~Demailly \cite{D10} -- \cite{D0}. The idea is to replace the
projective trace measure $T\wedge (dd^c\log|\cdot-x|)^p$ by
$T\wedge (dd^c\vp)^p$ with quite general psh functions $\vp$
(weights) with singularity at $x$. Classical and directional
Lelong numbers are particular cases of these ones, with specified
weight functions. Moreover, the technique of generalized Lelong
numbers allows to give more simple and natural proofs for deep
results concerning standard Lelong numbers.

The notion is based on machinery of complex Monge-Amp\`ere
operators, and eventually Monge-Amp\`ere currents are one of the
main objects in investigation of singularities of psh functions.
We start with a quick overview of basic results on Monge-Amp\`ere
operators, referring to \cite{B}, \cite{C}, \cite{D1}, \cite{Kol},
\cite{KoLN} for a more detailed introduction to the theory. The
presentation of generalized Lelong numbers follows Demailly's
paper \cite{D1}.

\subsection{Monge-Amp\`ere operators}

By {\it complex Monge-Amp\`ere operator} of a plurisubharmonic
function $u$ we mean the operator
$$
u\mapsto (dd^c u)^n
$$
or,
more generally, for plurisubharmonic functions $u_1,\ldots,u_p$,
$p\le n$, the wedge product
\beq
\label{eq:MAq}
dd^c u_1\wedge\ldots\wedge dd^c u_p.
\eeq
For $u$ smooth,
$$
(dd^cu)^n=\left({2\over \pi}\right)^n n!\det\left(
\frac{\partial^2u}{\partial z_j\partial\bar z_k}\right)\beta_n
$$
(and the corresponding mixed determinant for (\ref{eq:MAq})
with smooth $u_j$).

The problem is that such a wedge product cannot be defined for
arbitrary plurisubharmonic functions, see \cite{Siu2},
\cite{Kis7}. Using a Chern-Levine-Nirenberg estimate for the
complex Monge-Amp\`ere operator \cite{CLN},
 it was shown by Bedford and Taylor \cite{BT1}
that $(dd^cu)^k$ can be defined inductively as
$$
(dd^cu)^k=dd^c[u(dd^cu)^{k-1}]
$$
for continuous and even for all locally bounded psh functions $u$.
More generally, for any current $T\in\dpl$ and a function
$u\in PSH(\Omega)\cap L^\infty(\Omega)$, the current
$uT$ is well defined, has locally bounded mass, and
$dd^cu\wedge T:=dd^c(uT)\in {\cal D}_{p-1}^+$.

The following refined version of the Chern-Levine-Nirenberg estimate
is due to Cegrell.

\begin{prop}
\label{prop:CLN} {\rm \cite{C}} For any $K\Subset\Omega$ there is
a constant $C_{K,\Omega}$ such that for every $v\in PSH(\Omega)$
and $u_j\in PSH(\Omega)\cap L^\infty(\Omega)$, $1\le j\le p\le n$,
$$
\Vert vdd^c u_1\wedge\ldots\wedge dd^c u_p\Vert_K
\le C_{K,\Omega}\Vert v\Vert_{L^1(\Omega)}
\Vert u_1\Vert_{L^\infty(\Omega)}\ldots \Vert u_p\Vert_{L^\infty(\Omega)}.
$$
\end{prop}

So the obstacles for the definition of the operator arise from the
singularity sets of plurisubharmonic functions. For the
Monge-Amp\`ere operator be well defined, either the singularity
set has to be ``small'' or the function must not decrease too
rapidly to $-\infty$. The latter situation was studied in
\cite{B}, Theorem~4.3. However having in mind applications to
holomorphic mappings, one needs to make restrictions to the
singularity sets themselves (since the decay of $\log|f|$ is the
strongest possible). The following result is due to Fornaess and
Sibony.

\beth {\rm \cite{FoSi}} \label{theo:defMA} Let $T$ be a positive
closed current of bidimension $(p,p)$, and the unbounded loci
$L_j$ of plurisubharmonic functions $u_j$, $1\le j\le q\le p$,
satisfy \beq \label{eq:hausd} {\cal
H}_{2(p-m+1)}(L_{j_1}\cap\ldots\cap L_{j_m}\cap\supp T)=0 \eeq for
all choices of indices $j_1<\ldots<j_m$, $m=1,\ldots,q$, $ {\cal
H}_{2(p-m+1)}$ being the $2(p-m+1)$-dimensional Hausdorff measure.
Then the currents $u_1dd^cu_2\wedge\ldots\wedge dd^cu_q$ and
$$
dd^cu_1\wedge dd^cu_2\wedge\ldots\wedge dd^cu_q\wedge T:=
dd^c(u_1dd^cu_2\wedge\ldots\wedge dd^cu_q\wedge T)
$$
 are well defined and
have locally finite mass.
\eth

(The result with zero $(2p-2m+1)$-dimensional Hausdorff measure
was obtained earlier by Demailly \cite{D1}.)

Therefore, if $u_1,\ldots,u_q\in PSH(\Omega)\cap L^\infty_{loc}
(\Omega\setminus K)$, $K\Subset\Omega$, then $ dd^cu_1\wedge\ldots
dd^cu_k\wedge T\ge 0$ ($q\le p$) and in particular, $ (dd^cu)^q\in
{\cal D}_{n-k}^+(\Omega)$, is well defined for all $q\le n$. A
wider (though less explicitly presented) class of functions $u$
with well-defined Monge-Amp\`ere operator $(dd^cu)^n$ was
introduced in \cite{C1}, \cite{C2}.

Another problem is that the Monge-Amp\`ere operators are not
continuous with respect to the weak convergence of
plurisubharmonic functions (first mentioned in \cite{C0}). However
they are continuous under decreasing limits:

\beth {\rm \cite{D1}} \label{theo:contMA} Let functions $u_j$ and
a positive closed current $T$
 satisfy the conditions of Theorem~\ref{theo:defMA}.
If plurisubharmonic functions $u_j^s$ decrease to $u_j$ as $s\to\infty$,
then
$$
dd^cu_1^s\wedge dd^cu_2^s\wedge\ldots\wedge dd^cu_q^s\wedge T\to
dd^cu_1\wedge dd^cu_2\wedge\ldots\wedge dd^cu_q\wedge T
$$
in the weak topology of currents.
\eth


In particular, the result is true for functions locally bounded
outside a compact subset of $\Omega$.

When functions $u_j$ have the form $u_j=\log|f_j|$ with holomorphic
$f_j:\Omega\to\C$, condition (\ref{eq:hausd}) with $T=1$ means
\beq
\label{eq:genmtiple}
\dim  Z_{j_1}\cap\ldots\cap Z_{j_m}\le n-m,\quad m=1,2,\ldots,q,
\eeq
and for a function $u=\log\sum_j|f_j|^2=:\log|f|^2$,
the operator $(dd^cu)^q$ is well defined if
\beq
\label{eq:gendim}
\dim  Z_f\le n-q,
\eeq
 $Z_f=Z_1\cap\ldots\cap Z_q$ being the zero set of the mapping
$f=(f_1,\ldots,f_q)$.

For this specific situation, the following convergence result is known.

\beth \label{theo:convMAmap} {\rm \cite{Ro}, \cite{R1}, \cite{R2}}
Let a sequence of holomorphic mappings $f^s:\Omega\to\Cq$, $q\le
n$, converge to a mapping $f$ uniformly on compact subsets of
$\Omega$. If the zero set $Z_f$ of $f$ satisfies (\ref{eq:gendim})
then the currents $(dd^c\log|f^s|)^p$, $p\le q$, are well defined
on each subset $\omega\Subset\Omega$ for $s\ge s_0(\omega)$, and
$(dd^c\log|f^s|)^p\to (dd^c\log|f|)^p$,  $p\le q$. \eth

For locally bounded psh functions, the Monge-Amp\`ere operators
are continuous also for increasing sequences \cite{BT2}. More
general conditions can be given in terms of convergence with
respect to the Bedford-Taylor capacity \beq \label{eq:capacity}
Cap\,(E,\Omega):=\sup\{\int_E (dd^cu)^n:\: u\in PSH(\Omega),\:
-1\le u<0\}. \eeq A sequence of functions $u_s$ is said to
converge in capacity to $u$ on $\Omega$ if for any $\epsilon>0$
and any Borel set $K\Subset \Omega$,
$$
\lim_{s\to\infty}Cap\,(K\cap\{z:|u(z)-u_s(z)|>\epsilon\},\Omega)=0.
$$

\beth \label{theo:contcap} {\rm \cite{X1}} Let a sequence of
locally uniformly bounded plurisubharmonic functions
$\{u_j^s\}_{s=1}^\infty$ converge in capacity to $u_j$ as
$s\to\infty$, $1\le j\le q$, then
$$
dd^cu_1^s\wedge dd^cu_2^s\wedge\ldots\wedge dd^cu_q^s\to
dd^cu_1\wedge dd^cu_2\wedge\ldots\wedge dd^cu_q
$$
in the weak topology of currents.
\eth

The result remains true if $u_j^s$ are assumed to be uniformly
locally bounded outside a finite number of fixed points
$x_1,\ldots,x_m$ of $\Omega$, so they may have singularities at
these points \cite{X2}.

Some other sufficient conditions for the convergence are given in
\cite{FoSi}, \cite{C1}, \cite{C2}, \cite{X3}.

\subsection{Definition of generalized Lelong numbers}
\label{ssect:3.2}

Given $\varphi\in PSH(\Omega)$
and $r\in {\bf R}$, denote
$$
B_r(\varphi)=
\{z:\varphi(z)<r\},
$$
$$
S_r(\varphi)=\{z:\varphi(z)=r\}.
$$

A psh function $\varphi$ is {\it semiexhaustive} if
$B_R(\varphi)\Subset\Omega$ for some $R\in\R$.
In particular, $\varphi\in L^\infty_{loc}(\Omega\setminus
B_R(\varphi))$ and thus $ (dd^c\varphi)^k$ is well defined for all $k\le n$.
If, in addition, $\vp$ is such that $e^\vp$ is continuous on $\Omega$, it
is called a {\it psh weight} on $\Omega$.

Let $T\in\dpl$. Define
$$
\ntvr=\int_\brv T\wedge(dd^c\varphi)^p
$$
and
$$
\ntv=\lim_{r\to -\infty}\ntvr,
$$
{\it the generalized Lelong number}, or {\it the Lelong-Demailly number},
 with respect to the weight
$\varphi$ (\cite{D10}, \cite{D0}).

\medskip
\Exmp{s}.

1) $\varphi(z)=\log|z-x|\then\brv=B_{e^r}(x)$, $\ntvr=\nu(T,x,e^r)$ and
$\ntv=\ntx$.

2) the "directional" weights
\beq
\label{eq:dirw}
\varphi(z)=\varphi_{a,x}(z):=\sup_k a_k^{-1}\log|z_k-x_k|, \quad a_k>0,
\eeq
generate the directional Lelong
numbers with respect to $(a_1,\ldots,a_n)$ (to be shown in
Section~\ref{ssect:3.3}).

The following useful formula can be derived by means of Stokes' theorem.

\begin{prop}
\label{prop:formula} {\rm\cite{D1}} For any convex increasing
function $\gamma:\R\to\R$,
$$
\nu(T,\gamma\circ\vp,\gamma(r))=\gamma'(r)^p\ntvr,
$$
$\gamma'$ being understood as the left derivative.
In particular,
$$
\ntvr=e^{-2pr}\int_\brv T\wedge \left({1\over 2}
dd^ce^{2\varphi}\right)^p.
$$
\end{prop}

(Actually, we need not suppose $\gamma$ to be convex and
increasing, the only requirement being $\gamma\circ\vp$ to be a
psh weight \cite{Co}).

\subsection{Lelong-Jensen-Demailly formula}
\label{ssect:3.3}

 For the classical Lelong number of a psh function,
Theorem~\ref{theo:nu} serves as a bridge between its definition
as the density of the associated measure (i.e., as $\nu(dd^cu,x)$)
and as an asymptotic characteristic
of the function itself (as given in Corollary~\ref{cor:lognu}).
 A similar relation for the generalized Lelong numbers
exists, too.

Let $\vp$ be a psh weight in $\Omega$, $\vp_r=\max\,\{\vp,r\}$.
The swept out  Monge-Amp\`ere measure on $\srv$ is defined as
$$
\mrv=\left(dd^c\vp_r\right)^n|_{\srv}.
$$
It is a positive measure with the total mass
$\mrv(\srv)=(dd^c\vp)^n(\brv)$.
If $(dd^c\vp)^n=0$ on $\Omega\setminus\vp^{-1}(-\infty)$, then
$\mrv=(dd^c\vp_r)^n$.

\medskip
\Exmp. For $\vp=\log|z-x|$, $\mrv$ is the normalized Lebesgue measure on
$S_{e^r}(x)$.

More generally, if $\vp$ is smooth near $\srv$ and $d\vp\neq 0$ near
$\srv$, then
$$
\mrv=(dd^c\vp)^{n-1}\wedge d^c\vp|_{\srv}.
$$

\beth \label{theo:LJD} {\rm (Lelong-Jensen-Demailly formula)
\cite{D11}, \cite{D3}} Any $u\in PSH(\Omega)$ is
$\mu_r$-integrable for $-\infty<r<R$, and
$$
\mrv(u)-\int_\brv u(dd^c\vp)^n=\int_{-\infty}^r\nu(dd^cu,\vp,t)\,dt.
$$
\eth

\Exmp. For $\vp(z)=\log|z-x|$,  this becomes
$$
u(x)=
\lambda(u,x,\log r)
-\tau_{n-1}^{-1}\int_0^r{\sigma_u(x,t)\over t^{2n-1}}\,dt,\quad
0<r< {\rm dist}\,(x,\partial\Omega).
$$

As a consequence of the Lelong-Jensen-Demailly formula, we have

\beth
\label{theo:as}
 Let $(dd^c\vp)^n=0$ on
$\Omega\setminus\vp^{-1}(-\infty)$, then $r\mapsto\mrv(u)$ is a
convex function of $r$ and
$$
\nu(dd^c u,\vp)=\lim_{r\to
-\infty}{\mrv(u)\over r}.
$$
\eth

\Exmp{s}.

1) When $\vp(z)=\log|z-x|$, this is relation (\ref{eq:nu}) for the
classical Lelong numbers.

2) For $\vp(z)=\vp_{a,x}(z)$ defined by (\ref{eq:dirw}), we have
$$
B_r(\vp_{a,x})=\{z:|z_k-x_k|<e^{ra_k},\ \okn\},\quad
{\rm supp}\,\mrv=\{z:|z_k-x_k|=e^{ra_k},\ \okn\},
$$
and so $\mrv(u)=(a_1\ldots a_n)^{-1}\lambda(u,x,ra)$, which gives
us
$$
\nu(dd^cu,\vp_{a,x})=(a_1\ldots a_n)^{-1}\nu(u,x,a).
$$

\subsection{Semicontinuity properties}

The following results \cite{D1} are useful when studying families
of currents or weights.

\beth
\label{theo:scnt1}
If currents $T_k\in\dpl$ converge to a current $T$, then
$$
 \limsup_{k\to\infty}\nu(T_k,\vp)\le\ntv.
$$
\eth

\beth
\label{theo:scnt2}
If psh weights $\vp_k$ and $\vp$ are such that
$\exp\{\vp_k\}\to\exp\{\vp\}$ uniformly on compact subsets of
$\Omega$, then
$$
\limsup_{k\to\infty}\nu(T,\vp_k)\le\ntv.
$$
\eth

\subsection{Comparison theorems}

The first comparison theorem describes variation of the generalized Lelong numbers
with respect to the weights.
For a psh weight $\vp$, denote $L(\vp)=\vp^{-1}(-\infty)$.

\beth \label{theo:cmprs1} {\rm\cite{D3}, \cite{D1}}
 Let $T\in\dpl$ and
$\vp$ and $\psi$ be psh weights such that
$$
\limsup{\psi(z)\over\vp(z)}=l<\infty \quad {\rm as\ }{z\to L(\vp)},
\ z\in \supp T,
$$
then $\nu(T,\psi)\le l^p\ntv$.
\eth

The second comparison theorem indicates dependence of the generalized Lelong numbers
on the currents.

\beth \label{theo:cmprs2} {\rm\cite{D3}, \cite{D1}}
 Let $T\in\dpl$ and $u_k,\,v_k\in PSH(\Omega)$, $1\le
k\le q$, be such that the currents $dd^cu_1\wedge\ldots\wedge dd^cu_q
\wedge T$ and $dd^cv_1\wedge\ldots\wedge dd^cv_q\wedge T$ are well
defined (see Theorem~\ref{theo:defMA}),
$v_k=-\infty$ on $\supp T\cap L(\vp)$ and
$$
\limsup{u_k(z)\over v_k(z)}=l_k<\infty \quad
{\rm as\ }{z\to L(\vp)},
\ z\in \supp T\setminus v_k^{-1}(-\infty).
$$
Then $\nu(dd^cu_1\wedge\ldots\wedge dd^cu_q\wedge T,\vp)\le l_1\ldots
l_q\, \nu(dd^cv_1\wedge\ldots\wedge dd^cv_q\wedge T,\vp)$.
\eth

The above  results make it possible to obtain relatively simple
proofs for the statements of Theorems~\ref{theo:thie} --
\ref{theo:slice}. Besides, Theorem~\ref{theo:cmprs2} shows that
the residual Monge-Amp\`ere mass $(dd^cu)^n(x)$ of $u\in
PSH(\Omega)\cap L_{loc}^\infty(\Omega\setminus \{x\})$ is a
function of asymptotic behaviour of $u$ near $x$. This leads to
the notion of {\it standard singularities} at $x$ as equivalence
classes of psh functions with respect to their asymptotics
\cite{Za}. In particular, we have that the functions
$\max_j\log|z_j|^{a_j}$ and $\log\sum_j|z_j|^{a_j}$ (with
$a_1,\ldots,a_n>0$ being fixed) represent the same singularity at
$0$ and \beq\label{eq:ddclog} (dd^c\max_j\log|z_j|^{a_j})^n=
(dd^c\log\sum_j|z_j|^{a_j})^n=a_1\ldots a_n\, \delta_0. \eeq

Another application is the following result comparing the Lelong number of
a wedge product with the Lelong numbers of the factors.

\begin{cor}
\label{cor:estim} {\rm\cite{D1}} If $dd^cu_1\wedge\ldots\wedge
dd^cu_q$ is well defined, then
$$
\nu(dd^cu_1\wedge\ldots\wedge dd^cu_q\wedge T,x)\ge
\nu(u_1,x)\ldots\nu(u_q,x)\,\nu(T,x).
$$
\end{cor}

A remarkable relation between standard and generalized Lelong numbers
is given by

\begin{theo}
{\rm\cite{D1}} For $T\in\dpl$ and a psh weight $\vp$ with
$\vp^{-1}(-\infty)=x$,
$$
\ntp\ge\ntx\,\nu((dd^c\vp)^p,x).
$$
\end{theo}

When $p=1$ or $p=n-1$, this can be deduced from
the comparison theorems. But for $1<p<n-1$ they would give us a more rough
inequality
$\ntp\ge\ntx[\nu(\vp,x)]^p$. Instead, the proof in this case uses the relation
(assuming $x=0$)
$$
\ntp\ge\int_{U_n} \nu(T,\vp\circ g)\,d\kappa(g),
$$
$d\kappa$ being the Haar measure on the unitary group $U_n$, which
can be derived from the deep Theorem~\ref{theo:SiuDem}.


\section{Analyticity theorems for upperlevel sets}

Plurisubharmonicity assumes no {\it a priori} analyticity.
Nevertheless, analytic varieties appear from any psh function
(moreover, from any positive closed current) with singularities.

\subsection{Upperlevel sets for Lelong numbers}

Let $T\in\dpl$ and
$$
E_c(T):=\{x\in\Omega:\ntx\ge c\}, \ c>0,
$$
 be the upperlevel sets for the Lelong numbers of $T$.
Since $\ntx$ is lower semicontinuous, $E_c(T)$ is closed.
Furthermore, it
has locally finite ${\cal H}_{2p}$ Hausdorff
measure (this follows from Corollary~\ref{cor:volest}).

Section 4 is devoted to the following fundamental result:

\beth \label{theo:Siu} {\rm  Siu \cite{Siu}} $E_c(T)$ is an
analytic variety of dimension $\le p$. \eth

In other words, $\ntx$ is lower semicontinuous with respect to the
analytic Zariski topology (in which closed sets $=$ analytic varieties).

The original Siu's proof \cite{Siu} (1974) (developing results
from \cite{Bo} and \cite{Sk}) takes about 100 pages. A
considerable simplification was made by Lelong (1977) \cite{Le3}
who reduced the problem to that for psh function. In 1979 Kiselman
applied the attenuating singularities technique
(Subsection~\ref{ssect:1.6}) to get a more simple proof of Siu's
theorem for the classical  Lelong numbers \cite{Kis1} and, in
1986,  for the directional numbers \cite{Kis2}, \cite{Kis3}. His
ideas were used by Demailly to prove the theorem for the
generalized Lelong numbers \cite{D0} (1987). Perhaps, the shortest
known proof was proposed by Demailly in 1992 \cite{D2}.  It is
based on his approximation theorem for psh functions.

All the proofs rest heavily on $L^2$ estimates for the
$\bar\partial$ operator. We sketch below the proof of Theorem~\ref{theo:Siu}
based on Demailly's approximation theorem, as well as Kiselman-Demailly's proof
for generalized Lelong numbers.

\subsection{Reduction to plurisubharmonic functions}

The initial observation is as follows.
 Let $A=\{z: f_1(z)=\ldots=f_N(z)=0\}$, then the function
$v={1\over 2}\log\sum_k|f_k|^2$ has the property
$$
\nu([A],x)>0\iff x\in A\iff \nux>0.
$$
Actually, it is possible to construct a psh function whose Lelong numbers
just coincide with those of $A$ (and more generally, with the Lelong
numbers of an arbitrary given positive closed current).

\beth \label{theo:potent} {\rm \cite{Sk}, \cite{Sk2}, \cite{Le3}}
Let $\Omega$ be a pseudoconvex domain. Given a current $T\in\dpl$,
there exists a function
$u\in PSH(\Omega)$ such that $\nu(T,x)=\nux$ for every $x$.
\eth

The proof uses technique of canonical potentials
$$
U_j(z)=-\omega_p^{-1}\int|z-\zeta|^{-2p}\eta_j(\zeta)\,d\sigma_T(\zeta)
$$
with $\eta_j$ a non-negative, smooth function supported in $\Omega$,
and $\eta_j\equiv 1$ on a neighbourhood of $\overline\Omega_j\Subset\Omega$.
The function $U_j$ is subharmonic in ${\bf R}^{2n}$, and
$$
\sigma_{U_j}(x,r)={1\over 2\pi}\int_{B_r(x)}\Delta U_j=
[1+o(1)]\tau_{n-1}r^{2n-2}\ntx+ o(r^{2n-2}),\quad r\to 0,
$$
so
$$
\lim_{r\to 0}{\sigma_{U_j}(x,r)\over \tau_{n-1}r^{2n-2}}=\ntx \quad
\forall x\in\Omega_j.
$$
One can show that $dd^cU_j\ge -N_j \,dd^c|z|^2$, so
$u_j(z):=U_j(z)+N_j\,|z|^2+M_j\in PSH(\Omega)$ and $\nu(T,x)=\nu(u_j,x)$ for
all $x\in\Omega_j$. Exhausting $\Omega$ by $\Omega_j$ we get
the desired function $u$.

\subsection{$L^2$-extension theorems}

A bridge between plurisubharmonicity and analyticity is based on the
H\"ormander type results on solutions for the $\bar\partial$-problem.
In particular, the following two theorems have great importance in studying
singularities of plurisubharmonic functions.

\beth \label{theo:HBS} {\rm (H\"ormander-Bombieri-Skoda)
\cite{Bo}, \cite{Sk}, \cite{Sk1}} If $u$ is plurisubharmonic on a
pseudoconvex domain $\Omega$ and $e^{-u}\in L^2_{loc}(x)$ for some
$x\in\Omega$, then there exists a holomorphic function $f$ on
$\Omega$ such that
$$
\int |f|^2e^{-2u}(1+|x|^2)^{-m-\epsilon}\,\beta_p<\infty
$$
and $f(x)=1$.
\eth

\beth \label{theo:OT} {\rm (Ohsawa-Takegoshi) \cite{OT}} {Let $Y$
be an affine linear $p$-dimensional subspace of $\Cn$, $\Omega$ be
a bounded pseudoconvex domain in $\Cn$, and $u\in PSH(\Omega)$.
Then any function $h\in Hol(Y\cap\Omega)$ with
$$
\int_{Y\cap\Omega}|h|^2e^{-u}\,\beta_p<\infty
$$
can be extended to a function $f\in Hol(\Omega)$ and}
$$
\int_{\Omega}|f|^2e^{-u}\,\beta_n\le A(p,n,{\rm diam}\,\Omega)
\int_{Y\cap\Omega}|h|^2e^{-u}\,\beta_p.
$$
\eth

\subsection{Approximation theorem of Demailly}
\label{ssect:4.4}

Let $\Omega$ be a bounded pseudoconvex domain, $u\in PSH(\Omega)$.
Consider the Hilbert space
$$
H_m:=H_{m,u}(\Omega)=\{f\in Hol(\Omega):
\int_\Omega |f|^2e^{-2mu}\beta_n<\infty\},
$$
see Section~\ref{ssect:1.9}.
Let $\{\sigma_l^{(m)}\}_l$ be an orthonormal basis of $H_m$,
$$
u_m(z):={1\over 2m}\log\sum_l |\sigma_l^{(m)}(z)|^2\in PSH(\Omega).
$$
Note that $u_m(z)={1\over m}\sup\{\log|f(z)|:\,||f||_m<1\}$.

\beth \label{theo:appr} {\rm \cite{D2}, \cite{DK}} {There are
constants $C_1, C_2>0$ such that for any $z\in\Omega$ and every
$r<{\rm dist}\,(z,\partial \Omega)$,
$$
u(z)-{C_1\over m} \le u_m(z)\le \sup_{\zeta\in B_r(z)}u(\zeta)
+{1\over m}\log {C_2\over r^n}.
$$
In particular, $u_m\to u$ pointwise and in $L^1_{loc}(\Omega)$, and
\beq
\label{eq:appr}
\nux-{n\over m}\le\nu(u_m,x)\le\nux\quad\forall x\in\Omega.
\eeq
Besides, the integrability indices of $u$ and $u_m$ are related as
$$
I(u,x)-{1\over n}\le I(u_m,x)\le I(u,x).
$$
}
\eth

The proof uses Theorem~\ref{theo:OT} (more exactly, its particular case
of a one-point set $Y$).

It is worth mentioning that the functions $u_m$ from Theorem~\ref{theo:appr}
 control not only classical Lelong numbers of $u$ but also all its directional
ones \cite{R4} (see also \cite{FaJ1}): \beq \label{eq:apprdir}
\nuxa -m^{-1}\sum_j a_j\le \nu(u_m,x,a)\le\nuxa  \quad\forall
x\in\Omega,\ \forall a\in \Rn_+. \eeq

\subsection{Proof of Siu's theorem}
\label{ssect:prsiu}

By (\ref{eq:appr}),
$$
E_c(T)=E_c(u)=\bigcap_{m\ge m_0}E_{c-n/m}(u_m).
$$
Any set $E_a(u_m)$ is analytic since
$$
x\in E_{c-n/m}(u_m)\iff {\partial^\alpha\over\partial z^\alpha}
\sigma_l^{(m)}(x)=0\quad\forall\alpha:\: |\alpha|<cm-n,
$$
and so is $E_c(T)$.


\subsection{Siu's theorem for Lelong-Demailly numbers}

Due to relations (\ref{eq:apprdir}), Siu's theorem for directional
Lelong numbers can be proved exactly as for the classical ones in
Section~\ref{ssect:prsiu}. Another its proof was given by Kiselman
(see \cite{Kis3}) by means of the attenuating singularities
technique. It can be seen as a particular case of the following
Demailly's result on generalized Lelong numbers which actually is
a development of Kiselman's approach.

Let $X$ be a Stein manifold (e.g. a pseudoconvex domain in $\Cm$),
and $\vp$ be a semiexhaustive psh function on $\Omega\times X$.
The function $\vp_x(z):=\vp(z,x)$ is a psh weight on $\Omega$.
In this setting, $E_c=E_c(T,\vp)=\{x\in X:\nu(T,\vp_x)\ge c\}$.

\beth \label{theo:SiuDem} {\rm \cite{D1}}
 { If $\exp\vp\in C(\Omega\times X)$ and is locally H\"older
with respect to $x$, then $E_c$ is analytic in $X$.}
\eth

{\em Scheme of the proof.}

1) Construction of a family of psh potentials $u_a(x),\ a\ge 0$, whose
behaviour is determined by $\nu(T,\vp_x)$ (a refined version of
theorem~\ref{theo:potent}).

2) Let $N_{a,b}=\{x\in X:\exp\{-u_a/b\}\not\in L^2_{loc}(x)\}$,
then $E_c=\bigcap N_{a,b}$ where the intersection is taken
over all
$a<c$ and $b<(c-a)\gamma/m$ ($\gamma$ is the H\"older exponent of
$\vp$).

3) For any psh function $u$, the set
$$
NI_u=\{x\in X:e^{-u} \not\in L^2_{loc}(x)\}
$$
is analytic subset of $X$. This fact follows from Theorem~\ref{theo:HBS}.

4) Since $N_{a,b}=NI_{u_a/b}$, the conclusion follows.

\medskip
Partial Lelong numbers (\ref{eq:part})
 cannot be viewed as a particular case of Lelong-Demailly
numbers, and actually their upperlevel sets need not be analytic.
Some sufficient conditions for the analyticity were given in
\cite{W}, see also \cite{Kis3}.

\subsection{Semicontinuity theorems for integrability indices}
\label{ssect:4.7}

Let $I(u,x)$
be the integrability index of $u$ at $x$ (\ref{eq:intind}), and
$IE_c(u)=\{x:I(u,x)\ge c\}$ its upperlevel set.
Theorem~\ref{theo:HBS} implies Zariski's semicontinuity of the map
$x\mapsto I(u,x)$:

\beth
 $IE_c(u)$ is an analytic variety for all $c>0$.
\eth

The technique of Demailly's approximation theorem makes it possible
to show that the map $u\mapsto I(u,x)$ is upper-semicontinuous
with respect to the weak convergence of plurisubharmonic functions:

\beth {\rm \cite{DK}} If $\gamma>I(u,x)$ and $u_j\to u$, then
$\exp\{-u_j/\gamma\}\to \exp\{-u/\gamma\}$ in $L_{loc}^2$. \eth


\section{Structure of positive closed currents}

Importance of Siu's fundamental result (Theorem \ref{theo:Siu})
becomes clear by means of structural formulas for closed positive
currents. The proofs can be found in \cite{D1}.

\subsection{Siu's decomposition formula}

Let $T\in\dpl,\ A$ be an irreducible analytic variety of dimension $p$.
Define the value
$\nu(T,A):=\inf\,\{\ntx:x\in A\}$.
We have $\nu(T,A)=\ntx\ \forall x\in A\setminus A',\ A'$ being a proper
analytic subset of $A$. So, $\nu(T,A)$ is the {\it generic Lelong
number} of $T$ {\it along} $A$.

\begin{prop}
 $\chi_AT=\nu(T,A)\,[A]$ with $\chi_A$ the indicator
function of the set $A$.
\end{prop}

Note that  $\chi_AT\in\dpl$ in view of Theorem~\ref{theo:SkEM}.

\beth \label{theo:Siuform} {\rm (Siu's formula) \cite{Siu}} {For
any current $T\in\dpl$,
 there is a unique decomposition
$$
T=\sum_j\lambda_j[A_j]+R,
$$
where $\lambda_j>0$, $[A_j]$ are the integration currents over
irreducible $p$-dimensional varieties $A_j$, and $R\in\dpl$ is such
that $\dim E_c(R)<p$ for every $ c>0$.} Here  $\lambda_j$ are generic
Lelong numbers of $T$ along the varieties $A_j$.
\eth

Some components of lower dimension can actually occur in $E_c(R)$, however
 $\chi_{A}R=0$ for any $p$-dimensional variety $A$. And Siu's formula cannot
 be applied to $R$ since it is of bidimension $(p,p)$.

\subsection{King-Demailly formula}

Let $f=(f_1,\ldots,f_q):\Omega\to\Cq$ be a holomorphic mapping,
 $A_f=f^{-1}(0)$ be its zero set,
$$
u=\log|f|={1\over 2}\log\sum_k|f_k|^2.
$$
 If ${\rm dim}\,A_f
\le n-l$, then $(dd^cu)^l$ is well defined. Here
 ${\rm dim}\,A=\max_j\,{\rm dim}\,A_j$ and $A_j$ are
the irreducible components of $A_f$. Note also that $A_f=A_{f_1}\cap\ldots
\cap A_{f_q}$ where $A_{f_k}$ are the zero sets of the components $f_k$
of $f$ (of codimension $\le 1$).

When dealing with holomorphic mappings, it is convenient to
consider the corresponding {\it holomorphic chains}, i.e.,
the currents
$$
Z_f=\sum_j m_j[A_j]
$$
where the summation is made
 only for $(n-p)$-dimensional components $A_j$ of the variety
$A$ and $m_j\in{\bf Z}_+$ are the
generic multiplicities of $f$ at $A_j$.

The following result (originally established by J.~King for
$q=p\le n$) represents holomorphic chains as singular parts of
Monge-Amp\`ere currents.

\beth \label{theo:KD} {\rm \cite{GK}, \cite{D1}, \cite{Meo2}} If
the zero set $A_f$ of a holomorphic mapping $f:\Omega\to\Cq$ has
codimension $p$, then the currents $(dd^c\log|f|)^l$ and
$\log|f|(dd^c\log|f|)^l$ with $l<p$ have locally integrable
coefficients, and \beq \label{eq:KD} (dd^c\log|f|)^p=Z_f+R \eeq
  where $Z_f$ is the corresponding holomorphic chain and the current
 $ R\in{\cal D}_{n-p}^+(\Omega)$ is such that $\chi_{A_f}R=0$,
$\nu(R,x)$ are nonnegative integers  and $\{x:\nu(R,x)>0\}$ is an
analytic set of codimension at least $p+1$. \eth

{\em Particular cases}.

1) $p=q=1$: no condition on $A$ is required, $R=0$, so
we have the {\it Lelong-Poincar\`e equation}
$dd^c\log|f|=Z_f$.

2) $p=q\le n$: $R=0$, which gives {\it King's formula} \cite{GK}
$$
(dd^c\log|f|)^q=\sum_j m_j[A_j] =Z_f$$ (since $A$ has no
components of dimension $<n-q$).

\medskip

Another formula for the case $p=q\le n$ represents $Z_f$ as
the wedge product of the divisors of the
components $f_k$ of the mapping $f$. Denote $u_k=\log|f_k|$, then
$dd^cu_k=Z_{f_k}$.

\beth \label{theo:prod} {\rm\cite{D1}} If the zero sets $A_{f_k}$
satisfy condition (\ref{eq:genmtiple}), then
$$
dd^cu_1\wedge\ldots\wedge dd^cu_q= Z_{f_1}\wedge\ldots\wedge
Z_{f_q}=Z_f.
$$
\eth


\subsection{Lelong numbers of direct and inverse images }

Let $T\in\dpl$ and a holomorphic mapping $f:\Omega\to\Omega'$ be such that
its restriction to the support of $T$ is proper. Then its direct image
$f_*T\in{\cal D}_p^+(\Omega')$ and $\nu(f_*T,\vp)=\nu(T,\vp\circ f)$.
The problem is to evaluate this value in terms of characteristics of $T$ and $f$,
at least in the case of the classical Lelong numbers.

Let $x\in\Omega$ and $y=f(x)$. If the codimension of the fiber
$f^{-1}$ at $x$ is at least $p$, one can define the {\it multiplicity}
of $f$ at $x$ as
\beq
\label{eq:mup}
\mu_p(f,x):=\nu((dd^c\log|f-y|)^p,x).
\eeq
It is the multiplicity of the restriction of $f$ to a generic $p$-dimensional
plane through $x$ (and thus a non-negative integer).

\beth \label{theo:lbound} {\rm\cite{D1}} Let $T$ and $f$ be as
above. Let further $I(y)$ be the set of points $z\in\supp T\cap
f^{-1}(y)$ such that $z$ coincides with its connected component in
$\supp T\cap f^{-1}(y)$ and ${\rm codim}_z\,f^{-1}(y)\ge p$. Then
$$
\nu(f_*T,y)\ge\sum_{z\in I(y)}\mu_p(f,z)\nu(T,z).
$$
\eth

For finite holomorphic mappings, it is also possible to get an upper bound.
Set
$$
\bar\mu_p(f,z)=\inf_G\frac{[\sigma(G\circ f,z)]^p}{\mu_p(G,0)},
$$
where $G$ runs over all germs of maps $(\Omega,z)\to (\Cn,0)$ such that
$G\circ f$ is finite, and $\sigma(H,z)$ is the {\it Lojasiewicz
exponent} of a mapping $H$ at $z$
(the infimum of $\alpha>0$ such that $|\zeta-z|^\alpha/|H(\zeta)|$
is bounded near $z$).

\beth {\rm\cite{D1}} Let $f$ be a proper and finite holomorphic
mapping and $T\in\dpl$. Then
$$
\nu(f_*T,y)\le\sum_{z\in\supp T\cap f^{-1}(y)}\bar\mu_p(f,z)\nu(T,z).
$$
\eth

As an application, it gives the following relation for the Lelong
numbers with respect to the directional weights $\vp_{a,x}$
(\ref{eq:dirw}). Let $a_{j_1}\le\ldots \le a_{j_n}$, then
\beq\label{eq:dkbound} \frac{\ntx}{a_{j_{n-p+1}}\ldots
a_n}\le\nu(T,\vp_{a,x})\le \frac{\ntx}{a_{j_{1}}\ldots a_{j_p}}
\eeq
which is an extension of relation (\ref{eq:dircomp}) between
directional and classical Lelong numbers to currents of higher
degrees.

\medskip
Relations for the inverse images were studied in \cite{Meo}. Let
$f:\Omega\to\Omega'\subset\Cm$ be a surjective, proper and finite
holomorphic mapping. Then $\nu(f^*T,x)\ge \nu(T,y)$. For the
opposite direction, Theorem~\ref{theo:lbound} implies the bound
$$
\sum_{z\in f^{-1}(y)}\mu_p(f,z)\nu(f^*T,z)\le s\,\nu(T,y),
$$
where $s$ is the degree of the ramified covering $f$.

If $f$ is assumed to be only open and surjective, then
$$
\nu(f^*T,x)\le\mu_m(f,x)\nu(T,y).
$$

\subsection{Propagation of singularities}

Siu's theorem shows that once a psh function has logarithmic singularity
on a massive subset of an analytic variety, it must have them on the whole
variety. It gives an idea of constructing upper bounded psh functions
without (sub)extensions to larger domains.

\begin{prop}
{\rm\cite{RRo}} For every bounded convex domain $\Omega$ there
exists a negative function $u\in PSH(\Omega)$ with the following
property. If $\omega$ is an open subset of $\Omega$ and $v\in
PSH(\omega)$ satisfies $v\le u$ in $\omega$, then $v$ cannot be
extended to a function which is plurisubharmonic in a domain
strictly larger than $\Omega$.
\end{prop}
The obstacle for an extension here are upperlevel sets for the
Lelong numbers which are constructed to have no extension outside
the domain. On the other hand, for any negative psh function $f$
in the unit ball, any $\epsilon\in (0,1/n)$ and $r<1$ there exists
a psh function $v$ of logarithmic growth in $\Cn$ such that
$v(z)\le -|v(z)|^{{1\over n}-\epsilon}$ on $|z|<r$ \cite{EM1} (a
refined version of Josefson's theorem mentioned in
Section~\ref{ssect:psh}).

\medskip

Another application concerns unbounded maximal psh functions. A function
$u\in PSH(\Omega)$ is called {\it maximal} on $\Omega$ if for any
domain $\omega\Subset\Omega$ the relation $v\le u$ in
$\Omega\setminus\omega$ for $v\in PSH(\Omega)$ implies $v\le u$ in the
whole $\Omega$. For example, the function $\log|z_1|$ is maximal on
$\Cn$, $n>1$. Another example is $\log|f|$ with $f:\Omega\to\Cn$
such that its Jacobian determinant $J_f$ is identically zero. If
$J_f\not\equiv 0$,  $\log|f|$ is maximal outside the zero set $A_f$ of $f$.
Moreover, if the zero set contains components of positive dimension,
 $\log|f|$ can be maximal on a neighbourhood of some its points.

\beth {\rm\cite{R3}} For any holomorphic mapping $f:\Omega\to\Cn$
there exists a discrete set $CI_f\subset A_f$ such that  $\log|f|$
is maximal on a neighbourhood of every point $z\in\Omega\setminus
CI_f$. The set $CI_f$ is the complete indeterminence locus for $f$
considered as a meromorphic mapping to ${\bf P}^{n-1}$. \eth

So, $\log|f|$ is locally maximal on $\Omega\setminus CI_f$. It is not
known if it implies its maximality there (for bounded psh functions,
the maximality is a local property due to the characterization
$(dd^cu)^n=0$).

\section{Evaluation of residual Monge-Amp\`ere masses}

Here we will study the problem of evaluation of the Lelong numbers
of the Monge-Amp\`ere currents $(dd^cu)^m$. Even for $u=\log|f|$
with $f:\Omega\to\Cm$, $m>1$,  no explicit formulas for
computation of these values are available, and only bounds for
them are known (e.g., via the multiplicities of the components of
the mapping or in terms of the corresponding Newton polyhedra).

No upper bound of $ \nu((dd^cu)^m,x)$ in terms of $\nux$ is
possible. Nevertheless, it seems to be unknown if there exist a
psh function with zero Lelong number and nonzero residual
Monge-Amp\`ere mass.

As to lower estimates, a standard bound is
$$
 \nu((dd^cu)^m,x)\ge [\nux]^m
$$
 (cf. Corollary \ref{cor:estim}).
More precise relations can be obtained by means of more refined
characteristics of local behaviour of a function, e.g.,
directional Lelong numbers. To this end, we consider a notion of
local indicator, see \cite{LeR}.

\subsection{Definition and properties of local indicators}

Let $0\in\Omega$, $u\in PSH(\Omega), u\le 0$, and let $\nu(u,o,a)$
be its directional Lelong numbers (\ref{eq:dir}). Take
$t\in\Rn_-:=-\Rn_+$; the function
$$
\psi_u(t)=-\nu(u,0,-t)
$$ is
non-positive,
convex in $t$ and increasing in each $t_k$, so
$$
\Psi_u(z):=\psi_u(\log|z_1|,\ldots,\log|z_n|)
$$
is plurisubharmonic in $\{z:0<|z_k|<1,\: 1\le k\le n\}$ and thus
extends to a (unique) psh function in  the unit polydisk $D$,
which is called the {\it local indicator} of $u$ at $0$. (This
notion was introduced in \cite{LeR}, however plurisubharmonicity
of $\lambda(u,x,{\rm Re}z)$ was observed already in \cite{Kis1}.)

It is easily checked that
\beq\label{eq:hompsi}\psi_u(ct)=c\psi_u(t)\quad \forall c>0,\eeq
which implies $(dd^c\Psi_u)^n=0$ on $D\setminus\{z:z_1\cdot\ldots
z_n=0\}$.

Besides, $\Psi_{\Psi_u}=\Psi_u$, which means that $\Psi_u$ has the
same directional Lelong numbers as the function $u$.

\beth {\rm \cite{LeR}} For any function $u$ psh in a neighborhood
of $0$,
 \beq
\label{eq:ind_bound} u\le\Psi_u +C \eeq near the origin. \eth

\medskip
{\it Examples.}

1) For $u(z)=\log|z|$, $\Psi_u(z)=\sup_k\,\log|z_k|$.

2) If $u(z)=\vp_{a,0}(z)$ (\ref{eq:dirw}), then  $\Psi_u=\vp_{a,0}$,
so the directional weights are their own indicators.

3) Let $u=\log|f|,\ f:\Omega\to\Cm$, consider
the set
\beq
\label{eq:polyg}
\omega_0=\{J\in\Zn_+:\sum_j\left|{\partial^J  f_j\over\partial
z^J}(0)\right|\neq 0\}.
\eeq
As follows from (\ref{eq:index}),
 $\Psi_u(z)=\sup\,\{\log|z^J|: J\in\omega_0\}$.

\medskip

As was mentioned in Section \ref{ssect:tangents}, a psh function
need not have a unique tangent. At the same time, instead of the
family (\ref{eq:tang}) one may consider the collection
$$
 \tilde u_m(z)=m^{-1}u(z_1^m,\ldots,z_n^m),\quad m=1,2,\ldots.
$$

\beth {\rm \cite{R}} $\tilde u_m\to\Psi_u$ in $L_{loc}^1(D)$. \eth

So, local indicators can be viewed as {\it logarithmic tangents}
to psh functions.

\subsection{Reduction to indicators}

Let $\vp$ be a psh weight such that
$\vp^{-1}(-\infty)=0\in\Omega$.

In view of (\ref{eq:ind_bound}), Theorem \ref{theo:cmprs1} implies

\beth \label{theo:ind_bound} {\rm \cite{LeR}} If
$dd^cu_1\wedge\ldots\wedge dd^cu_q$ is well defined near the
origin (see Theorem~\ref{theo:defMA}), then
$$
\nu(dd^cu_1\wedge\ldots\wedge dd^cu_q,\vp)\ge
\nu(dd^c\Psi_{u_1}\wedge\ldots\wedge dd^c\Psi_{u_q},\vp)\ge
\nu(dd^c\Psi_{u_1}\wedge\ldots\wedge dd^c\Psi_{u_q},\Psi_\vp).
$$
\eth

For $u\in L^\infty_{loc}(\Omega\setminus\{0\})$,  the operator
$(dd^cu)^n$ is well defined, and the value
$$
\calr_u:=\nu((dd^cu)^n,0)
$$
is called the {\it residual measure} of $(dd^cu)^n$ at $0$. In
this situation, $(dd^c\Psi_u)^n=0$ on $D\setminus\{0\}$, so that
$$
(dd^c\Psi_u)^n=\tau_u\delta_0 $$
with
$$
\tau_u=\calr_{\Psi_u}
$$
which will be  called {\it the Newton number} of $u$ at $0$ (the
reason for using the name is clarified in the section
\ref{ssec:appl}).

\begin{cor}
\label{theo:nn} {\rm \cite{LeR}} If $u\in PSH(\Omega)\cap
L^\infty_{loc}(\Omega\setminus\{0\})$, then
 $\calr_u\ge \tau_u$.
\end{cor}

Similarly, for an $n$-tuple $(u)$ of psh function
$u_1,\ldots,u_n$, the residual Monge-Amp\`ere mass
$$\calr_{(u)}:=(dd^cu_1\wedge\ldots\wedge dd^cu_n)(0)$$
has the bound
$$ \calr_{(u)}\ge \tau_{(u)},$$
where
$$\tau_{(u)}=\calr_{(\Psi_u)}=(dd^c\Psi_{u_1}\wedge\ldots\wedge
dd^c\Psi_{u_n})(0).$$

To make all this reasonable, one has to look for good
bounds for the Newton numbers.

\subsection{Bounds in terms of directional Lelong numbers}

Let $\Psi(z)=\Psi(|z_1|,\ldots,|z_n|)\in PSH(D)\cap L_{loc}^\infty
(D\setminus 0)$, $\Psi< 0$ in $D$, and let its convex image
$\psi(t):=\Psi(\exp(t_1),\ldots,\exp(t_n))$ be homogeneous:
$\psi(ct)=c\psi(t),\ \forall c>0$. Such a function will be called
an {\it (abstract) indicator}. Note that our assumptions mean that
$\psi$ is the restriction to $\Rn_-$ of the support function of a
convex subset of $\Rn_+$. For a slightly different way to
introduce the abstract indicators, see \cite{Za0}, \cite{Za}.

We have $(dd^c\Psi)^n=\tau_\Psi\delta_0$. It is easy to see that
for all $z,\zeta\in D$, there is the inequality $\Psi(\zeta)\le
|\Psi(z)|\Phi_z(\zeta)$, where
$$
\Phi_z(\zeta)=\sup_k\,{\log|\zeta_k|\over |\log|z_k||}.
$$
Therefore, in view of Comparison Theorem~\ref{theo:cmprs2} and
equation (\ref{eq:ddclog}), we have for an $n$-tuple $(\Psi)$ of
indicators $\Psi_1,\ldots,\Psi_n$ the bound
$$
\calr_{(\Psi)}\ge |\Psi_1(z)\ldots\Psi_n(z)|(dd^c\Phi_z)^n=
\prod_k\left|\frac{\Psi_k(z)}{\log|z_k|}\right|
$$
for all $z\in D$ with $z_1\ldots z_n\neq 0$, and for a $q$-tuple
of indicators $\Psi_1,\ldots,\Psi_q$, $q<n$, and a psh weight
$\vp$,
$$
\nu(dd^cu_1\wedge\ldots\wedge dd^cu_q,\vp)\ge
|\Psi_1(z)\ldots\Psi_q(z)|\,\nu((dd^c\Phi_z)^q,\vp).
$$
Thus (\ref{eq:dkbound}) and Theorem \ref{theo:ind_bound} with
$\vp(z)=\log|z|$ imply the following result.

\beth \label{theo:dirbound} {\rm \cite{R}} In the conditions of
Theorem \ref{theo:ind_bound}, for any $q\le n$,
$$
\nu(dd^cu_1\wedge\ldots\wedge dd^cu_q,0)  \ge
\frac{\nu(u_1,0,a)\ldots\nu(u_q,0,a)}{a_{j_{q+1}}\ldots a_{j_n}}
\quad\forall a\in\Rn_+,
$$
where $a_{j_1}\le\ldots \le a_{j_n}$. In particular, if $u\in
PSH(\Omega)\cap L_{loc}^\infty(\Omega\setminus 0)$, then
$$
\calr_u\ge \tau_u\ge {[\nu(u,0,a)]^n\over a_1\ldots a_n}
\quad\forall a\in\Rn_+.
$$
\eth

\subsection{Geometric interpretation: volumes}

More sharp bounds can be obtained by precise calculation of the
Monge-Amp\`ere masses of the indicators.

Let $U(z)=U(|z_1|,\ldots,|z_n|)\in PSH(D)\cap L^\infty_{loc}
(D)$,
and $h(t):=U(\exp(t_1),\ldots,\exp(t_n))$ be its convex image
in $\Rn_-$. Then
$$
(dd^cU)^n=n!(2\pi)^{-n}{\cal MA}_\R[h]\,d\theta, \quad
(z_k=\exp\{t_j+i\theta_j\})
$$
${\cal MA}_\R$ being the {\it real Monge-Amp\`ere operator}. For
$h$ smooth,
$$
{\cal MA}_\R[h]=\det\left({\partial^2h\over\partial t_j\partial
t_k} \right)\,dt,
$$
and it extends,  as a measure-valued operator, to all convex
functions $h$. Furthermore, \beq \label{eq:rmavol} \int_F{\cal
MA}_\R[h]={\rm Vol}\,G_h (F),\quad F\Subset\Rn_-, \eeq with
$$
G_h(F)=\bigcup_{t^0\in F}\{a\in\Rn: h(t)\ge h(t^0)+
\langle a,t-t^0\rangle\ \forall t\in\Rn_-\}
$$
the {\it gradient image} of $F$ for the surface $\xi=h(t)$. The
terminology comes from the smooth situation because then we have
${\cal MA}_\R[h]$ equal the Jacobian determinant $J_{\nabla h}$
for the gradient mapping $\nabla h$, $G_h(F)= \nabla h(F)$ and the
equation (\ref{eq:rmavol}) is just the coordinate change formula.

 For a thorough treatment of the real
Monge-Amp\`ere operator, see \cite{RaT}.

So, for any $n$-circled set $E\Subset D$
and $E^* =\{t\in\Rn_-:
(\exp t_1,\ldots,\exp t_n)\in E\}$,
$$
(dd^cU)^n(E)=n!\,{\rm Vol}\,G_h(E^*).
$$
Let now $U=\max\{\Psi,-1\}$ with $\Psi$ an indicator. If $\Psi\in
L^\infty_{loc}(D\setminus\{0\})$, then the current $(dd^cU)^n$ is
supported by the set $E_\Psi=\{z:\Psi(z)=-1\}\Subset
D\setminus\{0\}$, and $E_\Psi^*=\Rn_+\setminus B_\Psi$, where
\beq\label{eq:bpsi}
B_\Psi=\{a\in\Rn_+:\langle
a,t\rangle\le\psi(t)\ \forall t\in\Rn_-\}
\eeq
(as before, $\psi$
is the convex image of $\Psi$). Note that the function $\psi$ is
just the restriction of the support function of the convex set
$B_\Psi$ to $\Rn_-$: $\psi(t)=\sup\,\{\langle a,t\rangle: a\in
B_\Psi\}$.

For any convex and complete subset $B$ of $\Rn_+$ (the latter
means that $a+\Rn_+\subset B$ for each $a\in B$) we put
$${\rm Covol}\,B= {\rm Vol}\,(\Rn_+\setminus B),
$$
the {\it covolume} of $B$.

\beth \label{theo:vol1} {\rm \cite{R} (see also \cite{D5})} If an
indicator $\Psi\in L^\infty_{loc}(D\setminus\{0\})$, then its
Monge-Amp\`ere mass is
$$\calr_\Psi = n!\,{\rm Covol}\,B_\Psi,
$$
where $B_\Psi$ is defined in (\ref{eq:bpsi}).
\eth

When $\Psi=\Psi_u$, the set
$$ B_\Psi=B_{u}= \{a\in\Rn_+:\nu(u,0,a)\le \langle a,b\rangle\ \forall
b\in\Rn_+\}.$$ Thus Corollary \ref{theo:nn} gives us

\beth\label{theo:vol2} If $u$ has isolated singularity at $0$,
then \beq\label{eq:vol} \calr_u\ge \tau_u= n!\,{\rm
Covol}\,B_{u}.\eeq \eth

Note that the value
$$
\sup_a\,{[\nu(u,0,a)]^n\over a_1\ldots a_n}
$$
from Theorem \ref{theo:dirbound} is the supremum over the volumes
of all simplices contained in the set $\Rn_+\setminus B_{u}$, and
$[\nu(u,0)]^n$ is the volume of the symmetric simplex
$\{a\in\Rn_+:\sum a_j\le\nu(u,0)\}\subset\Rn_+\setminus B_{u}$.

\medskip

To compute the mass of the corresponding mixed Monge-Amp\`ere
operators of indicators, we consider a (unique) form ${\rm
Covol}\,(B_1,\ldots,B_n)$ on $n$-tuples of complete convex subsets
$B_1,\ldots,B_n$ of $\Rn_+$ which is multilinear with respect to
Minkowsky's addition and such that for every convex complete $B$
with bounded complement in $\Rn_+$ we have ${\rm
Covol}\,(B,\ldots,B)={\rm Covol}\,B$. The form can be shown to be
well defined on all $n$-tuples $B_1,\ldots,B_n$ such that
$\Rn_+\setminus\cup_j B_j$ is bounded.

\beth\label{theo:vol3} If plurisubharmonic functions
$u_1,\ldots,u_n$ have $0$ as an isolated point of $\cap_k
\,L(u_n)$, then
$$(dd^cu_1\wedge\ldots\wedge dd^cu_n)(0) \ge n!\,{\rm
Covol}(B_{u_1},\ldots,B_{u_n}).$$
\eth

Taking $u_1=u$ and $u_2=\ldots=u_n=\varphi$, we get an estimate
for the Lelong-Demailly numbers $\nu(u,\vp)$. Another its form can
be derived by means of the Lelong-Jenson-Demailly formula (see
Theorem~\ref{theo:LJD}). Let $F$ be a subset of the convex set
$$L^\vp=\{t\in\Rn_-: \nu(\vp,0,-t)\ge 1\},$$
consider the sets
$$
\Gamma_F^\vp=\{a\in\Rn_+:\: \sup_{t\in F}\langle a,t\rangle =
\sup_{t\in L^\vp}\langle a,t\rangle=-1\},
$$
and
$$
\Theta_F^\vp=\{\gamma a:\: 0\le\gamma\le 1,\ a\in \Gamma_F^\vp\}.
$$

\beth \label{theo:mcw} {\rm\cite{R4}} If $\vp^{-1}(-\infty)=0$,
then \beq \label{eq:swept} \nu(u,\vp)\ge
n!\,\int_{E^\vp}\nu(u,0,-t)\,d\aoP(t), \eeq where the measure
$\aoP$ on the set $E^\vp$ of extreme points of $L^\vp$ is given by
the relation
$$
\aoP(F)={\rm Vol}\,\Theta_F^\vp
$$
for compact subsets $F$ of $E^\vp$.
\eth

\subsection{Functions with multicircled singularities}

We will say that a plurisubharmonic function $u$ on a domain
$\Omega\subset\Cn$ has {\it multicircled singularity} at a point
$x\in\Omega$ if there exists a multicircled plurisubharmonic
function $g$
 (i.e., $g(z)=g(|z_1|,\ldots,|z_n|)$) in a neighbourhood
of the origin such that \beq \exists\lim_{z\to x}{u(z)\over
g(z-x)}=1. \label{eq:amc} \eeq It is easy to see that $u$ has
multicircled singularity at $x$ if and only if it satisfies
relation (\ref{eq:amc}) with $g$ equal to some its
"circularization", say to
$$
g(z) = (2\pi)^{-n}
\int_{[0,2\pi]^n}u(x_1+z_1e^{i\theta_1},
\ldots,
x_n+z_ne^{i\theta_n})\,d\theta
$$
or to its maximum on the same set.

One can expect that such a regular behaviour implies nice relations
for the above characteristics of the singularity, and this is really the case.
As follows from Theorem~\ref{theo:intdir},
$$
I(u,x)=\sup\,\{\nuxa: \sum_ja_j=1\}
$$
 if $u$ has multicircled singularity at $x$.
Furthermore, the residual mass $\calr_u$ for such a function $u$
(assumed to be locally bounded outside the origin)  equals its
Newton number $\tau_u$ \cite{R5}, \cite{R4}. Moreover, it can be
estimated from above in terms of the limit values of the
directional Lelong numbers \cite{R5}: if
 $\nu'_j$ is the limit of
the values $\nu(u,0,a)$ as
$a_j\to 1$, $a_k\to +\infty$, $k\neq j$, then
$$
\calr_u=\tau_u\le \nu'_1\ldots\nu'_n.
$$
It implies a bound in terms of the corresponding partial Lelong
numbers $\nu_j$ (which can be strictly greater than $\nu'_j$ even
for multicircled functions, see Section \ref{ssec:pln}).

Finally, if a weight $\vp$ has multicircled singularity, then
there is also  equality in (\ref{eq:swept}), provided the
$-\infty$ set of $u$ does not contain lines parallel to the
coordinate axes and passing through $x$ \cite{R4}. And it fails to
be true, for example, for $u(z_1,z_2)=\log|z_1|$ and
$\vp(z_1,z_2)=\max\{-|\log|z_1||^{1/2},\log|z_2|\}+\log|z|$.

\subsection{Applications to holomorphic mappings: Newton
polyhedra}\label{ssec:appl}

For functions $u=\log|f|$ with holomorphic $ f:\Omega\to\Cn$ such
that $ f^{-1}(0)=0$, the residual Monge-Amp\`ere mass $\calr_u$ of
$(dd^cu)^n$ at $0$ is just the multiplicity $m_f$ of $f$ at the
origin.

Theorem \ref{theo:dirbound} with $a_1=\ldots=a_n$ gives us the
bound
$$
m_f\ge m_{f_1}\ldots m_{f_n}
$$
via the multiplicities of the components of the mapping $f$ (a
local variant of Bezout's theorem), and -- with integer $a_k$ --
the relation
$$
m_f\ge {m_{f_1,a}\ldots m_{f_n,a}\over a_1\ldots a_n}
$$
with the multiplicities of the weighted initial homogeneous
polynomial terms of the functions
$f_{j,a}(z)=f_j(z_1^{a_1},\ldots,z_n^{a_n})$ (Tsikh-Yuzhakov
\cite{YTs}, see also \cite{AYu}).

In view of (\ref{eq:index}), $B_{u}$ is the convex hull of the set
$\omega_0$ defined in (\ref{eq:polyg}), so it is the {\it Newton
polyhedron } for the mapping $f$ at $0$, and $\tau_u=N_f$ is the
{\it Newton number} of $f$, see \cite{AYu}, \cite{AGV},
\cite{Ku1}. This particular case of Theorem~\ref{theo:vol2}
recovers the result $m_f\ge N_f$ obtained by Kouchnirenko (1975)
by means of analytic and algebraic techniques.

The corresponding specification of Theorem~\ref{theo:vol3} gives a
Kouchnirenko-Bernstein's type result from \cite{AYu},
Theorem~22.10.

Note that for holomorphic mappings $ f:\Omega\to\Cq$ with ${\rm
codim}_0 f^{-1}(0)=q<n$, Theorem~\ref{theo:vol3} gives the bound
$$ m_f\ge n!\,{\rm Covol}(B_{u_1},\ldots,B_{u_q}, B_1,\ldots,
B_1),$$ where $u_j=\log|f_j|$, $B_{u_j}$ are the Newton polyhedra
of the functions $f_j$ at $0$, and $B_1=\{a\in\Rn_+: \sum_k a_k
\ge 1$\}.

Concerning Theorem~\ref{theo:mcw}, it is worth noticing that if a
weight $\vp$ has the form $\vp=\log|g|$ with a holomorphic mapping
$g$, then the measure $\aoP$ has finite support and is determined
by the $(n-1)$-dimensional faces of the Newton polyhedron of $g$.

So, methods of pluripotential theory are quite powerful  to
produce, in a simple and unified way, efficient bounds for
multiplicities of holomorphic mappings.

\vskip 0.5cm

Tek/Nat

H\o gskolen i Stavanger

POB 8002 Ullandhaug

4068 Stavanger

Norway

\vskip0.1cm

{\sc E-mail}: alexander.rashkovskii@tn.his.no

%
%
%
%
%
%
%
%
%
%

\end{document}